\documentclass[10pt]{article}
\newtheorem{teor}{Theorem}[section]
\newtheorem{lema}[teor]{Lemma}
\newtheorem{prop}[teor]{Proposition}

\newtheorem{defin}[teor]{Definition}
\newtheorem{obs}[teor]{Remark}

\usepackage[english]{babel}
\usepackage{amssymb,amsmath}
\usepackage{hyperref}

\title{The well adapted connection of a $(J^{2}=\pm 1)$-metric manifold}
\author{Fernando Etayo\footnote{Departamento de Matem\'{a}ticas, Estad\'{\i}stica y Computaci\'{o}n. Facultad de Ciencias.  Universidad de Cantabria. Avda. de los Castros, s/n, 39071 Santander, SPAIN. e-mail: etayof@unican.es}\, and Rafael Santamar\'{\i}a\footnote{Departamento de Matem\'{a}ticas. Escuela de Ingenier\'{\i}as Industrial e Inform\'{a}tica. Universidad de Le\'{o}n. Campus de Vegazana, 24071 Le\'{o}n, SPAIN. e-mail: rsans@unileon.es}}

\date{\today}
\begin{document}
\maketitle
\begin{abstract}
In this paper, we study the well adapted connection attached to  a $(J^{2}=\pm 1)$-metric manifold,  proving it exists for any of the four geometries and obtaining a explicit formula as a derivation law. Besides we  characterize the coincidence of the well adapted connection with the Levi Civita and the Chern connections.
\end{abstract}

{\bf Keywords:} $(J^{2}=\pm 1)$-metric manifold, natural connection, functorial connection, well adapted connection, Chern connection.

{\bf 2010 Mathematics Subject Classification:} 53C15, 53C50, 53C55, 53C07.

\vspace{7mm}

\section{Introduction}
\label{sec:introduction}

This paper concerns connections attached to $(J^{2}=\pm 1)$-metric manifolds. A manifold will be called to have an $(\alpha ,\varepsilon )$-structure if $J$ is an almost complex ($\alpha =-1$) or almost product ($\alpha =1$) structure  which is an isometry ($\varepsilon =1$) or anti-isometry ($\varepsilon =-1$) respect to $g$. Thus, there exist four kinds of $(\alpha ,\varepsilon )$ structures
\[
J^{2} = \alpha Id, \quad g(JX,JY)= \varepsilon g(X,Y), \quad \forall X, Y \in \mathfrak{X}(M). 
\]

As is well known, these four geometries have been intensively studied. The corresponding manifolds are known as:

\begin{enumerate}
\renewcommand*{\theenumi}{\roman{enumi})}
\renewcommand*{\labelenumi}{\theenumi}

\item  Almost-Hermitian manifold if it has a  $(-1,1)$-structure. We will consider through this paper the case $g$ being a Riemannian metric. See e.g.,  \cite{gray-hervella}.

\item Almost anti-Hermitian or almost Norden manifolds if it has a $(-1,-1)$-structure. The metric $g$ is semi-Riemannian having signature $(n,n)$. See e.g., \cite{ganchev-borisov}.

\item Almost product Riemannian manifolds if it has an $(1,1)$-structure. We will consider through this paper the case $g$ being a Riemannian metric, and the trace of $J$ vanishing, which in particular means these manifolds have even dimension. See e.g., \cite{naveira} and \cite{staikova}.

\item Almost para-Hermitian manifolds if it has an $(1,-1)$-structure. The metric $g$ is semi-Riemannian having signature $(n,n)$. See e.g., \cite{gadea}.
\end{enumerate}

\noindent The above cited papers focused on the classification of manifolds belonging to the different kind of geometries. Strong relations between these four geometries were established by us in \cite{debrecen}. Also, relations with
biparacomplex structures and 3-webs were shown, and we proved that any almost para-Hermitian structure admits an almost Hermitian structure with the same fundamental form.

If the structure $J$ is integrable, i.e., the Nijenhuis tensor $N_{J}=0$, the corresponding manifolds are called Hermitian, Norden, Product Riemannian and para-Hermitian (without the word ``almost"). Integrability means $M$ is a holomorphic manifold in cases i) and ii), and $M$ has two complementary foliations in cases iii) and iv).

As we have pointed out, the cases where $\varepsilon = 1$ the metric could be not Riemannian, as in  \cite{barros-romero} for the almost complex case and in \cite{pripoae} for the almost product one. We will focus only on the Riemannian metrics. Besides in all the cases, except the $(1,1)$ one, the trace of $J$ vanishes and the dimension of the manifold is even. We will impose in the case $(1,1)$ the condition $\mathrm{trace} \, J=0$ (of course, then $\mathrm{dim}\, M$ is also even) in order to have a unified treatment of these four geometries.

Observe that in the case $\alpha =1$  the tangent bundle of $M$
 can be decomposed (assuming $\mathrm{trace} \, J=0$) as the Whitney sum of two equidimensional subbundles corresponding to the eigenspaces of $J$;  i.e., $\mathrm{dim}\, T^+_p (M) = \mathrm{dim}\, T^-_p (M)$, for all $p \in M$, 
where
\[
T_p^+(M)=\{ v \in T_p(M) \colon Jv=v\}, \quad T_p^-(M)=\{ v \in T_p(M) \colon Jv=-v\}.
\]
The manifold $M$ is said  to have an almost paracomplex structure (see \cite{cruceanuetal}).   
\bigskip

In a lot of papers,  natural or adapted connections to an  $(\alpha ,\varepsilon )$-structure have been defined and studied. A \textit{natural} or \textit{adapted connection} is understood as a connection parallelizing both $J$ and $g$. Observe that the Levi Civita connection $\nabla^{\mathrm{g}}$ of $g$ is not a natural connection in the general case, although it is taken to define K\"{a}hler type condition as $\nabla^{\mathrm{g}} J =0$. On the other hand, as one can easily suppose, there is no a unique natural connection. In order to have a distinguished connection among all the natural ones, one must add some extra requirements. The first example of such a connection is that nowadays is called the Chern connection.
In 1946 Chern introduced $\nabla^{\mathrm{c}}$   as the unique linear connection in a Hermitian manifold $(M,J,g)$  such that  $\nabla^{\mathrm{c}} J=0$, $\nabla^{\mathrm{c}} g=0$ and having torsion $\mathrm{T}^{\mathrm{c}}$ satisfying
\[
\mathrm{T}^{\mathrm{c}}(JX,JY)= -\mathrm{T}^{\mathrm{c}}(X,Y), \quad \forall X, Y \in \mathfrak{X}(M).
\]

\noindent This connection also runs in the non integrable case, i.e., if $(M,J,g)$ is an almost Hermitian manifold, as it is showed; e.g., in \cite{gray}. Besides, Chern connection has been  taken as  a model in some other geometries, but possible ways to define connections depend on the values of $(\alpha ,\varepsilon)$.

\bigskip

In the present paper we will introduce a different approach to define  distinguished natural connections. We will focus the attention on the $G$-structure defined by each $(\alpha ,\varepsilon )$-structure, showing that all of them admit such a connection. This connection  $\nabla^{\mathrm{w}}$ is called \textit{the well adapted connection} and it is a functorial connection, see e.g., \cite{munoz2}.   The well adapted connection is in some sense the most natural connection one can define in a manifold with a $G$-structure, because it measures the integrability of the $G$-structure: it is integrable if and only if the torsion and the curvature tensors of the well adapted connection vanish (see \cite[Theor.\ 2.3]{valdes}). Surprisingly, the well adapted connection of an $(\alpha ,\varepsilon)$-structure has not been deeply studied. In \cite{valdes} and in \cite{brassov} well adapted connections to almost Hermitian (resp. almost para-Hermitian) manifolds were determined, but there are no results about the other two geometries. In this paper we  want to fill that lack  obtaining expressions of the well adapted connection to an  $(\alpha ,\varepsilon )$-structure in a unified way independent of the values  $(\alpha ,\varepsilon)$ as possible as we can do.

Thus, the main purposes of this paper are to prove the existence of the well adapted connection of any $(\alpha ,\varepsilon )$-structure, to obtain its expression as a derivation law, and to find the relationship between Chern type connections and the well adapted connection.  We will recover the results obtained by different mathematicians which will be named as the ``Bulgarian School'' through the present paper, because most of them are from that country. They have looked for distinguished natural connections in the case of almost Norden and almost product Riemannian manifolds  \textit{\`{a} la Chern}, i.e., imposing conditions on the torsion tensor field (see \cite{ganchev-mihova}, \cite{mihova}, \cite{staikova}, \cite{teofilova}). We will prove that their distinguished connections coincide with the well adapted connection in the corresponding geometry.

\bigskip

The organization of the paper is as follows:

In Section~\ref{sec:welladapted} we will recover the basic results about the well adapted connection, if there exists, of a $G$-structure. We will remember a sufficient condition for its existence, which will be the key to prove the existence of the well adapted connection of an $(\alpha ,\varepsilon)$-structure.

In Section~\ref{sec:existenceofwelladapted} we will obtain the $G$-structure associated to  $(J^{2}=\pm 1)$-metric-manifolds.  We will consider  $G_{(\alpha,\varepsilon )}$ as the $G$-structure associated to an  $(\alpha ,\varepsilon)$-structure $(J, g)$ and  we will prove the existence of the well adapted connection in each of the four $(\alpha ,\varepsilon)$ geometries (Theorem \ref{teor:ae-functorial}). The well adapted connection is found as a principal connection $\Gamma ^{\mathrm{w}}$ on the principal $G_{(\alpha,\varepsilon )}$-bundle.

In Section~\ref{sec:expressionofthewelladapted} we will obtain the expression of its derivation law $\nabla^{\mathrm{w}}$ (Theorem \ref{teor:bienadaptada-ae-estructura}). We will prove that $\nabla^{\mathrm{w}}$ is the unique adapted connection satisfying
\[
g(\mathrm{T}^{\mathrm{w}}(X,Y),Z)-g(\mathrm{T}^{\mathrm{w}}(Z,Y),X) = -\varepsilon (g(\mathrm{T}^{\mathrm{w}}(JX,Y),JZ)-g(\mathrm{T}^{\mathrm{w}}(JZ,Y),JX)),
\]

\noindent for all vector fields $X, Y, Z$ on $M$. Then, in Section~\ref{sec:particularizingthewelladaptedconnection}, we will specialize the expression in each of the four geometries, showing that the well adapted connections in almost Norden and almost product Riemannian manifolds coincide with connections previously introduced by the ``Bulgarian School" with different techniques. Finally we will characterize the equality between the well adapted and the Levi Civita connection (Theorem \ref{Levi Civita2}): they coincide if and only if $M$ is a K\"{a}hler-type manifold.

In Section~\ref{sec:chernconnection} we will show the classical Chern connection can be extended to the other geometry satisfying $\alpha \varepsilon = -1$, i.e., that of almost para-Hermitian manifolds, by the assumption (Theorem \ref{teor:chern-connection})
\[
\mathrm{T}^{\mathrm{c}}(JX,JY)= \alpha \mathrm{T}^{\mathrm{c}}(X,Y), \quad \forall X, Y \in \mathfrak{X}(M).
\]

\noindent Then the connection coincides with that introduced in \cite{etayo}. We will end the paper characterizing the coincidence between the Chern and the well adapted connection (Theorem \ref{teor:quasi}).

\bigskip

We will consider smooth manifolds and operators being of class $C^{\infty }$. As in this introduction, $\mathfrak{X}(M)$ denotes the module of vector fields of a manifold $M$. Endomorphisms of the tangent bundle are tensor fields of type (1,1). Almost complex and paracomplex structures are given by tensor fields of type (1,1). All the manifolds in the paper have even dimension $2n$. The general linear group will be denoted as usual by $GL(2n,\mathbb{R})$. The identity (resp. null) square matrix of order $n$ will be denoted as $I_n$ (resp. $O_n$).

Some proofs in the paper will seem to be reiterative, because one should check four similar cases. Nevertheless, these proofs need to be done with maximum care. In Remark \ref{noChern} we will show that the Chern-type connection which runs well in cases $\alpha \varepsilon =-1$, does not run in case $\alpha \varepsilon =1$. This is a warning to be extremely careful.

\section{The well adapted connection of a $G$-structure}
\label{sec:welladapted}

As is well known principal connections on the principal bundle of linear frames correspond to linear connections on the manifold, expressed by a derivation law. We will try to avoid the general constructions on frame bundles, but we need them to show what the well adapted connection means. We assume the theory of $G$-structures and reducible connections is known \cite{KN}.

The key results in order to define, if there exists, the well adapted connection of a $G$-structure are the following ones:

\begin{teor}[{\cite[Theor.\ 1.1]{valdes}}]
\label{teor:metodo}
Let $G\subseteq GL(n,\mathbb{R})$ be a Lie subgroup de Lie and let $\mathfrak g$ be its Lie algebra. The following two assertions are equivalent:

\begin{enumerate}
\renewcommand*{\theenumi}{\roman{enumi})}
\renewcommand*{\labelenumi}{\theenumi}

      \item For every $G$-structure $\pi \colon \mathcal S \to M$, there exists a unique linear connection
   $\Gamma^{\mathrm{w}}$ reducible to the $G$-structure such that, for every endomorphism  $S$ given by a section of the adjoint bundle of $\mathcal S$ and every vector field $X \in \mathfrak{X}(M)$ one has
            \[ \mathrm{trace}\, (S \circ i_X \circ \mathrm{T}^{\mathrm{w}}) = 0, \]
where $\mathrm{T}^{\mathrm{w}}$ is the torsion tensor of the derivation law  $\nabla ^{\mathrm{w}}$ of $\Gamma^{\mathrm{w}}$.

    \item If  $S \in \mathrm{Hom} (\mathbb{R}^{n} , \mathfrak g)$ satisfies
\[
    i_v \circ \mathrm{Alt} (S) \in \mathfrak g^{\perp},  \quad \forall v \in \mathbb{R}^{n},
    \]
   then $S = 0$, where  $\mathfrak g^{\perp}$ is the orthogonal subspace of $\mathfrak g$ in  $GL(n,\mathbb{R})$
    respect to the Killing-Cartan metric and
    \[ \mathrm{Alt} (S) (v,w) = S(v) w - S(w)v , \quad \forall v , w \in \mathbb{R}^{n}.\]
\end{enumerate}
The linear connection  $\Gamma^{\mathrm{w}}$, if there exists,  is called the well adapted connection to the $G$-structure $\pi \colon \mathcal S \to M$.
\end{teor}

This connection is a functorial connection in the sense of \cite{munoz2}. We do not develop the theory of functorial connection, looking for a direct introduction to the well adapted connection of an $(\alpha, \varepsilon )$-structure. Papers  \cite{munoz}, \cite{munoz2} and \cite{valdes}  cover that theory. In Section~\ref{sec:expressionofthewelladapted} we will explain and use all the elements introduced in the above Theorem.

\bigskip

The second result we need is following one:

\begin{teor}[{\cite[Theor.\ 2.4]{munoz2}, \cite[Theor. 2.1]{valdes}}]
\label{teor:suficiente}
 Let $G\subseteq GL(n,\mathbb{R})$ be a Lie group and let $\mathfrak g$ its Lie algebra. If $\mathfrak g^{(1)} = 0$ and $\mathfrak g$ is invariant under matrix transpositions, then condition   $ii)$ of Theorem \ref{teor:metodo} is satisfied, where $\mathfrak g^{(1)}=\{ S  \in \mathrm{Hom} (\mathbb{R}^{n}, \mathfrak g) \colon  S(v) w - S(w) v = 0,  \forall v, w \in \mathbb{R}^{n} \}$ denotes the first prolongation of the Lie algebra $\mathfrak g$.
\end{teor}

The well adapted connection, if there exists, is a functorial connection and measures the integrability of the $G$-structure, as it is proved in \cite[Theor.\ 2.3]{valdes}. The first condition also means that it is preserved by direct image. We do not use this property in the paper.

 \bigskip

 These results allow to check that there exists the well adapted connection to any $(\alpha ,\varepsilon )$-structure. In fact, the following one will be useful:

 \begin{prop}
 \label{teor:ppp-heridataria}
Let $G\subseteq GL(n,\mathbb{R})$ be a Lie group and let $\mathfrak g$ its Lie algebra such that $\mathfrak g^{(1)}=0$. If $H\subseteq G$ is a Lie subgroup then $\mathfrak h^{(1)}=0$, where $\mathfrak h$ denotes the Lie algebra of $H$.
\end{prop}

{\bf Proof.} As $H\subseteq G$ one also has $\mathfrak h \subseteq \mathfrak g$. Remember that
\[ \mathfrak h^{(1)}=
\{ S  \in \mathrm{Hom} (\mathbb{R}^{n}, \mathfrak h) \colon  S(v) w - S(w) v = 0,  \forall v, w \in \mathbb{R}^{n} \}.
\]
 Let $S \in \mathfrak h^{(1)}$. As $\mathfrak h \subseteq \mathfrak g$, then $S \in \mathrm{Hom} (\mathbb{R}^{n}, \mathfrak g)$, which implies
 \[
 S(v) w - S(w) v = 0,  \quad \forall v, w \in \mathbb{R}^{n},
 \]
 thus $S \in \mathfrak g^{(1)}=0$, which means $\mathfrak h^{(1)}=0$. $\blacksquare$

 \bigskip

In Section~\ref{sec:existenceofwelladapted} we  will show the $G$-structures corresponding to $(J^{2}=\pm 1)$-metric manifolds and in Section~\ref{sec:expressionofthewelladapted} we will prove they fulfill the sufficient conditions of Theorem \ref{teor:suficiente} to define the well adapted connection.

\section{Existence of the well adapted connection}
\label{sec:existenceofwelladapted}

We introduce the formal definition of $(\alpha, \varepsilon )$-structures.

\begin{defin}
Let $M$ be a  manifold, $J$ a  tensor field of type (1,1),  $g$ a semi-Riemannian metric in $M$ and $\alpha ,\varepsilon \in \{-1,1\}$. Then $(J, g)$ is called an $(\alpha ,\varepsilon )$-structure on $M$ if
\[
J^{2} = \alpha Id, \quad \mathrm{trace}\, J=0, \quad g(JX,JY)= \varepsilon g(X,Y), \quad \forall X, Y \in \mathfrak{X}(M),
\]

\noindent $g$ being a Riemanianan metric if $\varepsilon =1$. Then  $(M,J,g)$ is called a $(J^{2}=\pm 1)$-metric manifold.
\end{defin}

As we have said in the introduction, then $M$ has even dimension $2n$ and it is an almost complex (resp. almost paracomplex) manifold if $\alpha =-1$ (resp. $\alpha =1$).

\bigskip

We want to determine the $G$-structure associated to each $(\alpha, \varepsilon )$-structure. We use the following notation:

\begin{enumerate}

\item ${\mathcal F} (M)\to M$ is the principal bundle of linear frames of $M$; it has Lie structure  group $GL(2n,\mathbb{R})$.
\item ${\mathcal C}_{(\alpha ,\varepsilon )}\to M$ is the bundle of  adapted frames. It is a subbundle of the frame bundle.
\item $G_{(\alpha ,\varepsilon )}$ is the structure group of ${\mathcal C}_{(\alpha ,\varepsilon )}\to M$, i.e., the Lie group of the corresponding $G$-structure.
\item $\mathfrak{g}_{(\alpha, \varepsilon )}$ is the Lie algebra of $G_{(\alpha, \varepsilon )}$.
\end{enumerate}

\bigskip

Taking in mind the almost Hermitian case, which is well known, it is not difficult to find all the above elements. In fact, we have:

\begin{prop}[{\cite[Vol.\ II, Cap.\ IX]{KN}}]

Let $(J,g)$ be  a $(-1,1)$-structure in a $2n$-dimensional manifold $M$. The corresponding  $G_{(-1,1)}$-structure is given by
\[
{\mathcal C}_{(-1,1)} = \bigcup_{p \in M} \left\{ (X_1, \ldots, X_n, Y_1, \ldots , Y_n) \in {\mathcal F}_p (M) \colon
\begin{array}{l}
Y_i=JX_i,  \forall i =1, \ldots n, \\
g(X_i, X_j) = g(Y_i, Y_j)=\delta_{ij},\\
g(X_i, Y_j)=0, \forall i, j =1, \ldots, n
\end{array}
\right\},
\]

\[
G_{(-1,1)}=\left\{ \left( \begin{array}{cc}
A&B\\
-B&A
\end{array}\right) \in GL (2n; \mathbb{R}) \colon A, B \in GL (n; \mathbb{R}), \begin{array}{l}
A^t A +B^t B = I_n\\
B^t A-A^tB  = O_n
\end{array}  \right\},
\]
 and
\begin{equation}
\mathfrak{g}_{(-1,1)} = \left\{ \left( \begin{array}{cc}
A&B\\
-B&A
\end{array}\right) \in  \mathfrak{gl} (2n; \mathbb{R}) \colon  A, B \in  \mathfrak{gl} (n; \mathbb{R}),  \begin{array}{l}
A+A^t = O_n\\
B-B^t = O_n
\end{array}  \right\}.
\label{eq:(-11)-algebra}
\end{equation}
\end{prop}

\begin{prop}
Let $(J,g)$ be  a $(-1,-1)$-structure in a $2n$-dimensional manifold $M$. The corresponding  $G_{(-1,-1)}$-structure is given by
\[
{\mathcal C}_{(-1,-1)} = \bigcup_{p \in M} \left\{ (X_1, \ldots, X_n, Y_1, \ldots , Y_n) \in {\mathcal F}_p (M) \colon
\begin{array}{l}
Y_i=JX_i,  \forall i =1, \ldots n, \\
g(X_i, X_j) = g(Y_i, Y_j)=0,\\
g(X_i, Y_j)=\delta_{ij}, \forall i, j =1, \ldots, n
\end{array}
\right\},
\]
\[
G_{(-1,-1)}=\left\{ \left( \begin{array}{cc}
A&B\\
-B&A
\end{array}\right)\in GL (2n; \mathbb{R}) \colon  A, B \in GL (n; \mathbb{R}),  \begin{array}{l}
A^t A -B^t B = I_n\\
A^t B + B^t A  = O_n
\end{array}  \right\}
\]
and
\begin{equation}
\mathfrak{g}_{(-1,-1)} = \left\{ \left( \begin{array}{cc}
A&B\\
-B&A
\end{array}\right) \in  \mathfrak{gl} (2n; \mathbb{R}) \colon  A, B \in  \mathfrak{gl} (n; \mathbb{R}), \begin{array}{c}
A+A^t = O_n\\
B+B^t = O_n
\end{array}  \right\}.
\label{eq:(-1-1)-algebra}
\end{equation}
\end{prop}

\begin{prop}[{\cite{naveira}}]
Let $(J, g)$ be  a $(1,1)$-structure in a $2n$-dimensional manifold $M$. The corresponding  $G_{(1,1)}$-structure is given by
\[
{\mathcal C}_{(1,1)} = \bigcup_{p \in M}
\left\{ (X_1, \ldots, X_n, Y_1, \ldots , Y_n) \in {\mathcal F}_p (M) \colon
\begin{array}{l}
X_i \in T_p^+(M), Y_i \in T^-_p(M),  \forall i =1, \ldots n \\
g(X_i, X_j) = g(Y_i, Y_j)=\delta_{ij}\\
g(X_i, Y_j)=0, \forall i, j, =1, \ldots , n
\end{array}
 \right\},
\]

\[
G_{(1,1)}=\left\{ \left( \begin{array}{cc}
A&O_n\\
O_n&B
\end{array}\right) \in GL (2n; \mathbb{R}) \colon A, B \in O(n; \mathbb{R})    \right\},
\]
where $O(n; \mathbb{R})$ denotes the real orthogonal group of order $n$,
$O(n; \mathbb{R})=\{ N \in GL (n; \mathbb{R}) \colon N^t=N^{-1}\}$, and
\begin{equation}
\mathfrak{g}_{(1,1)}=\left\{ \left( \begin{array}{cc}
A&O_n\\
O_n&B
\end{array}\right) \in  \mathfrak{gl} (2n; \mathbb{R})  \colon A, B \in  \mathfrak{gl} (n; \mathbb{R})  \colon \begin{array}{l}
A  +A^t = O_n\\
 B + B^t  = O_n
\end{array}  \right\}.
\label{eq:(11)-algebra}
\end{equation}
\end{prop}

\begin{prop}[{\cite{cruceanuetal}, \cite{gadea}}]
Let $(J, g)$ be  a $(1,-1)$-structure in a $2n$-dimensional manifold $M$. The corresponding  $G_{(1,-1)}$-structure is given by
\[
{\mathcal C}_{(1,-1)} = \bigcup_{p \in M}
\left\{ (X_1, \ldots, X_n, Y_1, \ldots , Y_n) \in {\mathcal F}_p (M) \colon
\begin{array}{l}
X_i \in T_p^+(M), Y_i \in T^-_p(M),  \forall i =1, \ldots n \\
g(X_i, X_j) = g(Y_i, Y_j)=0\\
g(X_i, Y_j)=\delta_{ij}, \forall i, j, =1, \ldots , n
\end{array}
 \right\},
\]

\[
G_{(1,-1)}=\left\{ \left( \begin{array}{cc}
A&O_n\\
O_n&(A^t)^{-1}
\end{array}\right) \in GL (2n; \mathbb{R}) \colon A \in GL (n; \mathbb{R}) \right\},
\]
and
\begin{equation}
\mathfrak{g}_{(1,-1)}=\left\{ \left( \begin{array}{cc}
A&O_n\\
O_n&-A^t
\end{array}\right) \in  \mathfrak{gl} (2n; \mathbb{R})  \colon A \in  \mathfrak{gl} (n; \mathbb{R})  \right\}.
\label{eq:(1-1)-algebra}
\end{equation}
\end{prop}

In the four cases, we have emphasized the Lie  algebras, because they are the key elements to prove the existence of the well adapted connection. Now we are going to study the relationships among different structural groups we have found. The same relationships will have their Lie algebras. Let us remember the orthogonal and neutral orthogonal Lie groups and algebras:
\begin{eqnarray*}
O(2n; \mathbb{R})&=&\{N \in GL (2n; \mathbb{R}) \colon  N^t = N^{-1}\},\\
\mathfrak{o} (2n; \mathbb{R}) &=& \{ N \in  \mathfrak{gl} (2n; \mathbb{R}) \colon N+N^t=O_{2n}\},\\
O(n,n;\mathbb{R})&=&\left\{
\left(
\begin{array}{cc}
A&B\\
C&D
\end{array}
\right)
\in GL (2n; \mathbb{R}) \colon A, B, C, D \in GL (n; \mathbb{R}),
\begin{array}{l}
C^tA+A^tC=O_n\\
C^tB+A^tD=I_n\\
D^tB+B^tD=O_n
\end{array}
\right\},\\
\mathfrak{o}(n,n;\mathbb{R})&=&\left\{
\left(
\begin{array}{cc}
A&B\\
C&D
\end{array}
\right)
\in  \mathfrak{gl} (2n; \mathbb{R}) \colon A, B, C, D \in GL (n; \mathbb{R}),
\begin{array}{l}
A+D^t=O_n\\
B+B^t=O_n\\
C+C^t=O_n
\end{array}
\right\}.
\end{eqnarray*}

The following result will be relevant in order to prove the existence of the well adapted connection to an $(\alpha ,\varepsilon )$-structure.

\begin{prop}
\label{teor:prolongacionesnulas}
Let $n \in \mathbb N$. The first prolongation of the Lie algebras of $O(2n;\mathbb{R})$ and $O(n,n;\mathbb{R})$ vanish; {\em i.e.};
\[
\mathfrak o(2n;\mathbb{R})^{(1)}=0, \quad  \mathfrak o(n,n;\mathbb{R})^{(1)}=0.
\]
\end{prop}

 {\bf Proof.} One can prove the result by a straightforward computation. In a different way, one can deduce the result from some properties of functorial connections:
\begin{enumerate}

\item The well adapted connection is a functorial connection.

\item Manifolds endowed with a Riemannian or a neutral metric admits the well adapted connection, which is the Levi Civita connection (see \cite[Theor.\ 3.1 and Rem.\ 3.2]{munoz}).

\item If a manifold admits a functorial connection associated to a $G$-structure, then $\mathfrak g^{(1)}=0$ (see  \cite[Theor.\ 2.1]{munoz2}).

\end{enumerate}
The  result is proved. $\blacksquare$

 \bigskip

One can easily prove the following relationships:

\begin{prop}
\label{contenidos}
 Assuming the above notations, one has the following subsets:

\begin{enumerate}

\item $G_{(-1,1)}\subseteq O(2n; \mathbb{R}), \quad G_{(1,1)}\subseteq O(2n;\mathbb{R}), \quad G_{(-1,-1)}\subseteq O(n,n;\mathbb{R}), \quad  G_{(1,-1)} \subseteq O(n,n;\mathbb{R}).$

\item $\mathfrak{g}_{(-1,1)}\subseteq \mathfrak{o}(2n; \mathbb{R}), \quad \mathfrak{g}_{(1,1)}\subseteq \mathfrak{o}(2n;\mathbb{R}), \quad
\mathfrak{g}_{(-1,-1)}\subseteq \mathfrak{o}(n,n;\mathbb{R}) \quad  \mathfrak{g}_{(1,-1)} \subseteq \mathfrak{o}(n,n;\mathbb{R}).$

\end{enumerate}

Besides, one has the following equalities:
\begin{eqnarray*}
G_{(-1,1)}= GL (n; \mathbb C ) \cap O(2n; \mathbb{R}), &\quad& G_{(1,1)}=GL (n; \mathbb{R})\times GL (n; \mathbb{R}) \cap O(2n;\mathbb{R}),\\
G_{(-1,-1)}=GL (n; \mathbb C ) \cap O(n,n;\mathbb{R}),&\quad& G_{(1,-1)} = GL (n; \mathbb{R})\times GL (n; \mathbb{R}) \cap O(n,n;\mathbb{R}).
\end{eqnarray*}

\end{prop}

Taking into account the above results we can prove:

\begin{teor}
\label{teor:ae-functorial}
Let $M$ be a $2n$-dimensional manifold with an $(\alpha ,\varepsilon )$-structure. Then $M$ admits the well adapted connection.
\end{teor}

 {\bf Proof.} Let $(\alpha ,\varepsilon ) \in \{(-1,1), (-1,-1), (1,1,), (1,-1)\}$. The four Lie algebras  $\mathfrak g_{(\alpha ,\varepsilon )}$ corresponding to the four Lie groups  $G_{(\alpha ,\varepsilon)}$ are given by formulas (\ref{eq:(-11)-algebra}), (\ref{eq:(-1-1)-algebra}), (\ref{eq:(11)-algebra}) y (\ref{eq:(1-1)-algebra}). They are invariant under matrix transpositions, as one can easily check. Taking into account Proposition \ref{contenidos}, we obtain, by
 Proposition \ref{teor:ppp-heridataria} and  Proposition \ref{teor:prolongacionesnulas}, the formula $\mathfrak g_{(\alpha,\varepsilon)}^{(1)}=0$  holds. Then, by  Theorem \ref{teor:suficiente} we conclude the result. $\blacksquare$
\bigskip

In Section~\ref{sec:expressionofthewelladapted} we will study carefully this well adapted connection. In order to pass from the bundles to the manifolds, we need the following result, similar to that of the  case  given in \cite[Vol. II, Prop. 4.7]{KN}.

\begin{prop}
\label{teor:aeconexiones}
Let $(M, J, g)$ be a  $(J^2=\pm1)$-metric manifold and let $\pi \colon {\mathcal C}_{(\alpha ,\varepsilon )} \to M$ the  $G_{(\alpha ,\varepsilon )}$-structure on $M$  defined by $(g,J)$. Let $\Gamma$ a linear connection on  $M$ and let $\nabla$ the corresponding derivation law. Then $\Gamma$ is a reducible connection to $\pi \colon {\mathcal C}_{(\alpha ,\varepsilon )} \to M$ if and only if $\nabla J = 0, \nabla g = 0$.
\end{prop}

Thus, reducible connections correspond to natural or adapted connections. Among them, we have a distinguished one: the well adapted. In the next Section we will study the well adapted connection as a derivation law.

\bigskip

\begin{obs}
The above results could suggest that every known $G$-structure admits a well adapted connection. This is not the case. For example, the $G$-structure determined by an almost complex or an almost paracomplex structure does not admit a well adapted connection. In both cases one can use an idea taken in the proof of Proposition \ref{teor:prolongacionesnulas}: if there exists a functorial connection associated to a $G$-structure, then $\mathfrak g^{(1)}=0$. The corresponding structural groups are $GL(n,\mathbb{C})$ and $GL(n,\mathbb{R})\times GL(n,\mathbb{R})$. The first prolongation of the corresponding Lie algebras does not vanish, thus proving these $G$-structures do not admit a functorial connection; in particular, the well adapted connection.
\end{obs}

\section{Expression of the  well adapted connection}
\label{sec:expressionofthewelladapted}

Let $(M, J, g)$ be a $(J^{2}=\pm 1)$-metric manifold and let $\nabla^{\mathrm{w}}$ be the covariant operator or derivation law defined by the well adapted connection of the corresponding $G_{(\alpha ,\varepsilon )}$-structure. Taking into account Proposition \ref{teor:aeconexiones} we know $\nabla^{\mathrm{w}} J=0$ and  $\nabla^{\mathrm{w}} g=0$. As $\nabla^{\mathrm{w}}$ is uniquely determined, there should exist another condition which determine it uniquely. We are looking for this condition. So, the natural way is to study the set of all the natural connections of a $(J^{2}=\pm 1)$-metric manifold.

\begin{defin}
Let $(M, J, g)$ be a $(J^2=\pm1)$-metric manifold.  A  covariant derivative or derivation law $\nabla^{\mathrm{a}}$ on $M$ is said to be natural or adapted to $(J, g)$ if
\[
\nabla^{\mathrm{a}} J=0, \quad  \nabla^{\mathrm{a}} g=0;
\]
i.e., if it is the derivation law of a linear connection on $M$ reducible to the $G_{(\alpha ,\varepsilon )}$-structure on $M$ defined by $(J, g)$.
\end{defin}

\begin{defin}
\label{teor:tensorpotencial}
Let $(M,J,g)$ be a $(J^2=\pm1)$-metric manifold, let  $\nabla^{\mathrm{g}}$ be the derivation law of the Levi Civita connection of $g$ and let $\nabla^{\mathrm{a}}$ be a derivation law adapted to $(J,g)$. The potential tensor of $\nabla^{\mathrm{a}}$ is the  tensor $S\in \mathcal T^1_2 (M)$ defined as
\[
S(X,Y)=\nabla^{\mathrm{a}}_X Y -\nabla^{\mathrm{g}}_X Y, \quad \forall X, Y \in \mathfrak{X} (M).
\]
\end{defin}

Then we can parametrize all the adapted connections to $(J,g)$:

\begin{lema}
\label{teor:natural}
Let $(M,J,g)$ be a $(J^2=\pm1)$-metric manifold. The set of derivation laws adapted to $(J,g)$ is:
\[
\left\{ \nabla^{\mathrm{g}}+S \colon S \in \mathcal T^1_2 (M),
\begin{array}{l}
JS(X,Y)-S(X,JY)=(\nabla^{\mathrm{g}}_X J) Y, \\
 g(S(X,Y),Z)+g(S(X,Z),Y)=0,
 \end{array}
 \quad
\forall X, Y, Z \in \mathfrak{X} (M) \right\}.
\]
\end{lema}

{\bf Proof.} Let $S$ be the potential tensor of  $\nabla^{\mathrm{a}}$; then
\[
\nabla^{\mathrm{a}}_X Y = \nabla^{\mathrm{g}}_X Y + S(X,Y), \quad \forall X, Y \in \mathfrak{X} (M).
\]
Given $X, Y \in \mathfrak{X} (M)$, if the following relations hold
\begin{eqnarray*}
\nabla^{\mathrm{a}}_X J Y &=& J \nabla^{\mathrm{a}}_X Y \\
\nabla^{\mathrm{g}}_X JY + S(X,JY) &=& J \nabla^{\mathrm{g}}_X Y + J S(X,Y)
\end{eqnarray*}
then the condition  $\nabla^{\mathrm{a}} J=0$ is equivalent to the following identity
\[
JS(X,Y)-S(X,JY) =\nabla^{\mathrm{g}}_X JY - J\nabla^{\mathrm{g}}_X Y =(\nabla^{\mathrm{g}}_X J) Y, \quad \forall X, Y \in \mathfrak{X} (M).
\]
Given the vector fields  $X, Y, Z$ on $M$, as $\nabla^{\mathrm{g}} g=0$ one has
\begin{eqnarray*}
(\nabla^{\mathrm{a}}_X g) (Y,Z) &=& (\nabla^{\mathrm{a}}_X g) (Y,Z)- (\nabla^{\mathrm{g}}_X g) (Y,Z)\\
                         &=& -g(\nabla^{\mathrm{a}}_X Y, Z)-g(\nabla^{\mathrm{a}}_X Z, Y)+ g(\nabla^{\mathrm{g}}_X Y, Z)+g(\nabla^{\mathrm{g}}_X Z, Y)\\
                         &=&-(g(S(X,Y),Z)+g(S(X,Z),Y).
\end{eqnarray*}
thus proving  $\nabla^{\mathrm{a}}$ parallelizes the metric $g$ if and only if
\[
g(S(X,Y),Z)+g(S(X,Z),Y)=0, \quad \forall X, Y, Z \in \mathfrak{X} (M). \ \blacksquare
\]

The fundamental result in this section is the following one:

 \begin{teor}
 \label{teor:bienadaptada-ae-estructura}
 Let $(M,J,g)$ be a $(J^2=\pm1)$-metric manifold. The  derivation law $\nabla^{\mathrm{w}}$ of the well adapted connection $\Gamma ^{\mathrm{w}}$ is the unique derivation law satisfying $\nabla^{\mathrm{w}} J=0, \nabla^{\mathrm{w}} g=0$ and
\begin{equation}
g(\mathrm{T}^{\mathrm{w}}(X,Y),Z)-g(\mathrm{T}^{\mathrm{w}}(Z,Y),X) = -\varepsilon (g(\mathrm{T}^{\mathrm{w}}(JX,Y),JZ)-g(\mathrm{T}^{\mathrm{w}}(JZ,Y),JX)), \quad \forall X,Y, Z \in \mathfrak{X} (M).
\label{eq:welladapted}
\end{equation}
\end{teor}

The relations $\nabla^{\mathrm{w}} J=0$, $\nabla^{\mathrm{w}} g=0$ hold because of  Proposition \ref{teor:aeconexiones}, being $\Gamma ^{\mathrm{w}}$ reducible to the corresponding $G_{(\alpha, \varepsilon)}$-structure. The hard part of the proof is formula (\ref{eq:welladapted}). The strategy we will follow is this: we must prove that condition $i)$ in Theorem \ref{teor:metodo} is equivalent to formula (\ref{eq:welladapted}), because that condition characterizes the well adapted connection. This can be done working with local adapted frames defined in local charts (which will be introduced in Definition \ref{teor:ae-baselocal})  but first we should show what a section of the adjoint bundle means in our context of  $G_{(\alpha ,\varepsilon )}$-structures. This is our first auxilar result.

\bigskip

Let $\pi \colon \mathcal C_{(\alpha,\varepsilon)}\to M$ be the bundle of adapted frames, and let  $\mathrm{ad} \mathcal C_{(\alpha,\varepsilon)} =(\mathcal C_{(\alpha,\varepsilon)} \times \mathfrak g_{(\alpha,\varepsilon)})/G_{(\alpha,\varepsilon)}$ be the adjoint bundle. The structural group   $G_{(\alpha,\varepsilon)}$ acts on $\mathcal C_{(\alpha,\varepsilon)} \times \mathfrak g_{(\alpha,\varepsilon)}$ as:

\[
(u_p, A) \cdot N = (u_p \cdot N, N^{-1} A N), \quad \forall p\in M, u_p \in {\mathcal C_{(\alpha,\varepsilon)}} _p,  N \in G_{(\alpha,\varepsilon)}, A \in \mathfrak g_{(\alpha,\varepsilon)}.
\]
As $\mathfrak g_{(\alpha,\varepsilon)}\subseteq  \mathfrak{gl} (n; \mathbb{R})$, one has
$\mathrm{ad} \mathcal C_{(\alpha,\varepsilon)} \subseteq \mathrm{End} (TM)$, where $\mathrm{End} (TM)$ denotes  the set of endomorphisms of the tangent bundle of $M$. Then we have:

\begin{prop}
\label{teor:adjunto-ae-estructura}
Let $(M, J, g)$ be a $(J^2=\pm1)$-metric manifold of dimension $2n$. The sections of the adjoint bundle $\mathrm{ad} \mathcal C_{(\alpha,\varepsilon)}$ are the endomorphisms of the tangent bundle of $M$ satisfying the following two conditions:
\[
J\circ S = S \circ J, \quad g(SX, Y)=-g(X,SY), \quad \forall X, Y \in \mathfrak{X} (M).
\]
\end{prop}

{\bf Proof.} Observe that given   $p\in M$, then an element $S \in (\mathrm{ad} \mathcal C_{(\alpha,\varepsilon)})_p$ is an endomorphism $S\colon T_p (M)\to T_p(M)$ having coordinate matrix belonging to the Lie algebra $\mathfrak g_{(\alpha,\varepsilon)}$ when it is expressed  respect to the adapted frame  $u_p=(X_1, \ldots, X_n, Y_1, \ldots, Y_n) \in {\mathcal C_{(\alpha,\varepsilon)}}_p$. Then we will prove both implications working at a point $p\in M$.\medskip

$\Rightarrow)$ We will first consider the case $\alpha=1$ with the two subcases $\varepsilon =\pm 1$ and  after that we will take $\alpha=-1$ with the corresponding subcases. The proof follows from a careful analysis of the four situations. Condition on $\alpha$ will determine the commutativity of $J$ and $S$. Condition on $\varepsilon$ will allow us to obtain the formula linking $g$ and $S$.
\medskip

Assuming $\alpha=1$ there exist  $A, B \in  \mathfrak{gl} (n; \mathbb{R})$ such that the endomorphism $S$ has the following matrix respect to the adapted frame  $u_p$
\[
\left(
\begin{array}{cc}
A&O_n\\
O_n &B
\end{array}
\right) \in  \mathfrak{gl} (2n; \mathbb{R}),
\]
i.e.,
\[
SX_j = \sum_{i=1}^n a_{ij}X_i, \quad SY_j = \sum_{i=1}^n b_{ij} Y_i, \quad \forall j=1, \ldots , n,
\; {\rm where}\;
JX_i =X_i, \quad JY_i =-Y_i, \quad \forall i=1, \ldots, n.
\]
Then, for each $j =1, \ldots, n$, one has
\begin{eqnarray*}
JSX_j = J \left( \sum_{i=1}^n a_{ij}X_i\right) =  \sum_{i=1}^n a_{ij}X_i, &\quad&
JSY_j = J \left( \sum_{i=1}^n b_{ij}Y_i\right) = - \sum_{i=1}^n b_{ij}Y_i,\\
SJX_j= SX_j = \sum_{i=1}^n a_{ij}X_i, &\quad&
SJY_j= -SY_j =  -\sum_{i=1}^n b_{ij} Y_i,
\end{eqnarray*}
thus proving $(J\circ S) (u) = (S \circ J) (u)$ for all $u\in T_p(M)$.
\medskip

If $\varepsilon =1$ then $A+A^t=O_n$ and $B+B^t=O_n$, thus obtaining
\[
a_{ij}=-a_{ji}, \quad b_{ij}=-b_{ji}, \quad \forall i, j =1,\ldots , n.
\]
where
\[
g(X_i,X_j)=g(Y_i, Y_j)= \delta_{ij}, \quad g(X_i, Y_j)= 0, \quad \forall i, j =1, \ldots, n.
\]
For each pair $i, j =1, \ldots , n$, one has
\begin{eqnarray*}
g(SX_i, X_j)= g\left(\sum_{j=1}^n a_{ji}X_j, X_j \right)= a_{ji}, &\quad & g(X_i, SX_j)= g\left(X_i, \sum_{i=1}^n a_{ij}X_i \right)= a_{ij}=-a_{ji},\\
g(SY_i, Y_j)= g\left(\sum_{j=1}^n b_{ji}Y_j, Y_j \right)= b_{ji}, &\quad& g(Y_i, SY_j)= g\left(Y_i, \sum_{i=1}^n b_{ij}Y_i \right)= b_{ij}=-b_{ji}\\
g(SX_i, Y_j)= g\left(\sum_{j=1}^n a_{ji}X_j, Y_j \right)= 0, & \quad & g(X_i, SY_j)= g\left(X_i, \sum_{i=1}^n b_{ij}Y_i \right)= 0,
\end{eqnarray*}
thus proving $g(Su, v)=-g(u, Sv)$ for all  $u, v \in T_p (M)$.
\medskip

If $\varepsilon =-1$ then $B=-A^t$, thus obtaining
\[
b_{ij}=-a_{ji}, \quad \forall i, j =1,\ldots , n.
\]
where
\[
g(X_i,X_j)=g(Y_i, Y_j)= 0, \quad g(X_i, Y_j)= \delta_{ij}, \quad \forall i, j =1, \ldots, n.
\]
For each pair $i, j =1, \ldots , n$, one has
\begin{eqnarray*}
g(SX_i, X_j)= g\left(\sum_{j=1}^n a_{ji}X_j, X_j \right)= 0, &\quad& g(X_i, SX_j)= g\left(X_i, \sum_{i=1}^n a_{ij}X_i \right)= 0,\\
g(SY_i, Y_j)= g\left(\sum_{j=1}^n b_{ji}Y_j, Y_j \right)= 0, &\quad& g(Y_i, SY_j)= g\left(Y_i, \sum_{i=1}^n b_{ij}Y_i \right)= 0,\\
g(SX_i, Y_j)= g\left(\sum_{j=1}^n a_{ji}X_j, Y_j \right)= a_{ji}, &\quad& g(X_i, SY_j)= g\left(X_i, \sum_{i=1}^n b_{ij}Y_i \right)= b_{ij}=-a_{ji},
\end{eqnarray*}
thus proving $g(Su, v)=-g(u, Sv)$ for all  $u, v \in T_p (M)$.
\medskip

Assuming $\alpha=-1$ there exist $A, B \in  \mathfrak{gl} (n; \mathbb{R})$ such that the endomorphism $S$ has the following matrix respect to the adapted frame  $u_p$
\[
\left(
\begin{array}{cc}
A&B\\
-B &A
\end{array}
\right) \in  \mathfrak{gl} (2n; \mathbb{R}),
\]
i.e.,
\[
SX_j = \sum_{i=1}^n a_{ij}X_i -\sum_{i=1}^n b_{ij} Y_i, \quad SY_j = \sum_{i=1}^n b_{ij} X_i+ \sum_{i=1}^n a_{ij} Y_i, \quad \forall j=1, \ldots , n,
\]
where
\[
JX_i =Y_i, \quad JY_i =-X_i, \quad \forall i=1, \ldots, n.
\]
For each $j =1, \ldots, n$, one has
\begin{eqnarray*}
JSX_j &=& J \left( \sum_{i=1}^n a_{ij}X_i -\sum_{i=1}^n b_{ij} Y_i \right) =  \sum_{i=1}^n b_{ij}X_i +\sum_{i=1}^n a_{ij} Y_i , \\
JSY_j &=& J \left( \sum_{i=1}^n b_{ij} X_i+ \sum_{i=1}^n a_{ij} Y_i\right) = -\sum_{i=1}^n a_{ij} X_i+ \sum_{i=1}^n b_{ij} Y_i,\\
SJX_j&=& SY_j = \sum_{i=1}^n b_{ij} X_i+ \sum_{i=1}^n a_{ij} Y_i, \\
SJY_j&=& -SX_j =  -\sum_{i=1}^n a_{ij}X_i +\sum_{i=1}^n b_{ij} Y_i,
\end{eqnarray*}
thus proving $(J\circ S) (u) = (S \circ J) (u)$ for all $u\in T_p(M)$.
\medskip

If $\varepsilon =1$ then $A+A^t=O_n$ and $B-B^t=O_n$, thus obtaining
\[
a_{ij}=-a_{ji}, \quad b_{ij}=b_{ji}, \quad \forall i, j =1,\ldots , n.
\]
where
\[
g(X_i,X_j)=g(Y_i, Y_j)= \delta_{ij}, \quad g(X_i, Y_j)= 0, \quad \forall i, j =1, \ldots, n.
\]
For each pair $i, j =1, \ldots , n$, one has
\begin{eqnarray*}
g(SX_i, X_j)= g\left(\sum_{i=1}^n a_{ji}X_j -\sum_{i=1}^n b_{ji} Y_j, X_j \right)= a_{ji},
&\quad &
g(X_i, SX_j)= g\left(X_i, \sum_{i=1}^n a_{ij}X_i -\sum_{i=1}^n b_{ij} Y_i \right)= a_{ij}=-a_{ji},\\
g(SY_i, Y_j)= g\left(\sum_{i=1}^n b_{ji}X_j +\sum_{i=1}^n a_{ji} Y_j, Y_j \right)= a_{ji},
&\quad&
g(Y_i, SY_j)= g\left(Y_i, \sum_{i=1}^n b_{ij} X_i+ \sum_{i=1}^n a_{ij} Y_i \right)= a_{ij}=-a_{ji},\\
g(SX_i, Y_j)= g\left(\sum_{i=1}^n a_{ji}X_j -\sum_{i=1}^n b_{ji} Y_j, Y_j \right)= -b_{ji},
&\quad&
g(X_i, SY_j)= g\left(X_i,\sum_{i=1}^n b_{ij} X_i+ \sum_{i=1}^n a_{ij} Y_i \right)= b_{ij}=b_{ji},
\end{eqnarray*}
thus proving $g(Su, v)=-g(u, Sv)$ for all $u, v \in T_p (M)$.
\medskip

If $\varepsilon =-1$ then $A+A^t=O_n$ and $B+B^t=O_n$, thus obtaining
\[
a_{ij}=-a_{ji}, \quad b_{ij}=-b_{ji}, \quad \forall i, j =1,\ldots , n.
\]
where
\[
g(X_i,X_j)=g(Y_i, Y_j)= 0, \quad g(X_i, Y_j)= \delta_{ij}, \quad \forall i, j =1, \ldots, n.
\]
For each pair $i, j =1, \ldots , n$, one has
\begin{eqnarray*}
g(SX_i, X_j)= g\left(\sum_{i=1}^n a_{ji}X_j -\sum_{i=1}^n b_{ji} Y_j, X_j \right)= -b_{ji},
&\quad&
g(X_i, SX_j)= g\left(X_i, \sum_{i=1}^n a_{ij}X_i -\sum_{i=1}^n b_{ij} Y_i \right)= -b_{ij}=b_{ji},\\
g(SY_i, Y_j)= g\left(\sum_{i=1}^n b_{ji}X_j +\sum_{i=1}^n a_{ji} Y_j, Y_j \right)= b_{ji},
&\quad&
g(Y_i, SY_j)= g\left(Y_i, \sum_{i=1}^n b_{ij} X_i+ \sum_{i=1}^n a_{ij} Y_i \right)= b_{ij}=-b_{ji},\\
g(SX_i, Y_j)= g\left(\sum_{i=1}^n a_{ji}X_j -\sum_{i=1}^n b_{ji} Y_j, Y_j \right)= a_{ji},
&\quad&
g(X_i, SY_j)= g\left(X_i,\sum_{i=1}^n b_{ij} X_i+ \sum_{i=1}^n a_{ij} Y_i \right)= a_{ij}=-a_{ji},
\end{eqnarray*}
thus proving $g(Su, v)=-g(u, Sv)$ for all $u, v \in T_p (M)$.
\medskip

$\Leftarrow)$ Let  $S$ be an endomorphism of the tangent bundle of $M$. Its matrix respect to an adapted frame $u_{p}\in {\mathcal C_{(\alpha,\varepsilon)}}_p$ will be

\[
\left(
\begin{array}{cc}
A&B\\
C &D
\end{array}
\right) \in  \mathfrak{gl} (2n; \mathbb{R}),
\quad {\rm i.e.,}\quad
SX_j=\sum_{i=1}^n a_{ij} X_i + \sum_{i=1}^n c_{ij} Y_i, \quad SY_j = \sum_{i=1}^n b_{ij} X_i + \sum_{i=1}^n d_{ij} Y_i, \quad \forall j=1, \ldots, n.
\]

We must prove that if $S$ satisfies the two relations with $J$ and $g$, then $S$ is a section of the adjoint bundle, which is equivalent to prove that the above matrix of $S$ belongs to the corresponding Lie algebra $\mathfrak g_{(\alpha,\varepsilon)}$. As in the other implication we begin assuming $\alpha=1$ with the two subcases $\varepsilon =\pm 1$ and after that we will take $\alpha =-1$ with the corresponding subcases.
\medskip

Let $\alpha=1$. Then $JX_i =X_i, JY_i =-Y_i, \forall i=1, \ldots, n$, thus obtaining
\begin{eqnarray*}
SJX_j= \sum_{i=1}^n a_{ij} X_i + \sum_{i=1}^n c_{ij} Y_i,
&\quad&
SJY_j=-\sum_{i=1}^n b_{ij} X_i - \sum_{i=1}^n d_{ij} Y_i, \\
JSX_j= \sum_{i=1}^n a_{ij} X_i - \sum_{i=1}^n c_{ij} Y_i,
&\quad&
JSY_j= \sum_{i=1}^n b_{ij} X_i - \sum_{i=1}^n d_{ij} Y_i.
\end{eqnarray*}

As $J\circ S = S \circ J$ one has $c_{ij}=0, d_{ij}=0, \forall i, j =1\ldots, n$,
thus proving the matrix of $S$ has the following expression:
\[
\left(
\begin{array}{cc}
A&O_n\\
O_n &D
\end{array}
\right) \in  \mathfrak{gl} (2n; \mathbb{R}),
\quad {\rm i.e.}, \quad
SX_j=\sum_{i=1}^n a_{ij} X_i  \quad SY_j =  \sum_{i=1}^n d_{ij} Y_i, \quad \forall j=1, \ldots, n.
\]

Let $\varepsilon=1$. Thus $g(X_i,X_j)=g(Y_i, Y_j)= \delta_{ij}, g(X_i, Y_j)= 0, \forall i, j =1, \ldots, n$, and then for each pair $i, j =1, \ldots , n$, one has
\begin{eqnarray*}
g(SX_i, X_j)= g \left( \sum_{j=1}^n a_{ji} X_j , X_j \right)=a_{ji},
&\quad&
g(X_i, SX_j)= g \left( X_i, \sum_{j=1}^n a_{ij} X_i \right)=a_{ij},\\
g(SY_i, Y_j)= g \left( \sum_{j=1}^n d_{ji} Y_j , Y_j \right)=d_{ji},
&\quad&
g(Y_i, SY_j)= g \left( Y_i, \sum_{j=1}^n d_{ij} Y_i \right)=d_{ij}.
\end{eqnarray*}
As $g(Su, v)=-g(u, Sv)$ for all $u, v \in T_p (M)$ one has
$a_{ij}=-a_{ji}, \quad d_{ij}=-d_{ji}, \quad \forall i, j=1, \ldots, n$, thus proving the matrix of $S$ belongs to
 $\mathfrak g_{(1,1)}$ (see (\ref{eq:(11)-algebra})).
 \medskip

Let $\varepsilon=-1$. Thus $g(X_i,X_j)=g(Y_i, Y_j)= 0, g(X_i, Y_j)= \delta_{ij}, \forall i, j =1, \ldots, n$, and then for each pair $i, j =1, \ldots , n$, one has
\begin{eqnarray*}
g(SX_i, Y_j)= g \left( \sum_{j=1}^n a_{ji} X_j , Y_j \right)=a_{ji}, \quad
g(X_i, SY_j)= g \left( X_i, \sum_{j=1}^n d_{ij} Y_i \right)=d_{ij}.
\end{eqnarray*}
As $g(Su, v)=-g(u, Sv)$ for all $u, v \in T_p (M)$ one has
$a_{ij}=-d_{ji},  \quad \forall i, j=1, \ldots, n$, thus proving the matrix of $S$ belongs to $\mathfrak g_{(1,-1)}$ (see (\ref{eq:(1-1)-algebra})).
\medskip

Let $\alpha=-1$. Then $JX_i =Y_i,  JY_i =-X_i, \forall i=1, \ldots, n$, thus obtaining
\begin{eqnarray*}
SJX_j \sum_{i=1}^n b_{ij} X_i + \sum_{i=1}^n d_{ij} Y_i,
&\quad&
SJY_j=-\sum_{i=1}^n a_{ij} X_i - \sum_{i=1}^n c_{ij} Y_i,\\
JSX_j= -\sum_{i=1}^n c_{ij} X_i a +\sum_{i=1}^n a_{ij} Y_i,
&\quad&
JSY_j=-\sum_{i=1}^n d_{ij} X_i + \sum_{i=1}^n b_{ij} Y_i.
\end{eqnarray*}

As $J\circ S = S \circ J$ one has $b_{ij}=-c_{ij}, d_{ij}=a_{ij}, \forall i, j =1\ldots, n$, thus proving the matrix of $S$ has the following expression:
\[
\left(
\begin{array}{cc}
A&B\\
-B &A
\end{array}
\right) \in  \mathfrak{gl} (2n; \mathbb{R}),
\quad {\rm i.e.,}\quad
SX_j=\sum_{i=1}^n a_{ij} X_i - \sum_{i=1}^n b_{ij} Y_i, \quad SY_j = \sum_{i=1}^n b_{ij} X_i + \sum_{i=1}^n a_{ij} Y_i, \quad \forall j=1, \ldots, n.
\]

Let $\varepsilon=1$. Thus
$g(X_i,X_j)=g(Y_i, Y_j)= \delta_{ij}, g(X_i, Y_j)= 0, \forall i, j =1, \ldots, n$, and then for each pair
 $i, j =1, \ldots , n$, one has
\begin{eqnarray*}
g(SX_i, X_j)= g \left( \sum_{j=1}^n a_{ji} X_j -\sum_{j=1}^n b_{ji} X_j , Y_j \right)=a_{ji},
&\quad&
g(X_i, SX_j)= g \left( X_i, \sum_{j=1}^n a_{ij} X_i - \sum_{i=1}^n b_{ij} Y_i \right)=a_{ij},
\\
g(SX_i, Y_j)= g \left( \sum_{j=1}^n a_{ji} X_j -\sum_{j=1}^n b_{ji} Y_j, Y_j \right)=-b_{ji},
&\quad&
g(X_i, SY_j)= g \left( X_i, \sum_{j=1}^n b_{ij} X_i + \sum_{i=1}^n a_{ij} Y_i  \right)=b_{ij}.
\end{eqnarray*}
As $g(Su, v)=-g(u, Sv)$ for all $u, v \in T_p (M)$ one has $a_{ij}=-a_{ji}, b_{ij}=b_{ji}, \forall i, j=1, \ldots, n$, thus proving the matrix of $S$ belongs to  $\mathfrak g_{(-1,1)}$ (see (\ref{eq:(-11)-algebra})).
 \medskip

Finally, let  $\varepsilon=-1$. Thus $g(X_i,X_j)=g(Y_i, Y_j)= 0,  g(X_i, Y_j)= \delta_{ij},  \forall i, j =1, \ldots, n$, and then for each pair  $i, j =1, \ldots , n$, one has
\begin{eqnarray*}
g(SX_i, X_j)= g \left( \sum_{j=1}^n a_{ji} X_j -\sum_{j=1}^n b_{ji} X_j , Y_j \right)=-b_{ji},
&\quad&
g(X_i, SX_j)= g \left( X_i, \sum_{j=1}^n a_{ij} X_i - \sum_{i=1}^n b_{ij} Y_i \right)=-b_{ij},
\\
g(SX_i, Y_j)= g \left( \sum_{j=1}^n a_{ji} X_j -\sum_{j=1}^n b_{ji} Y_j, Y_j \right)=a_{ji},
&\quad&
g(X_i, SY_j)= g \left( X_i, \sum_{j=1}^n b_{ij} X_i + \sum_{i=1}^n a_{ij} Y_i  \right)=a_{ij}.
\end{eqnarray*}
As $g(Su, v)=-g(u, Sv)$, for all $u, v \in T_p (M)$ one has $
a_{ij}=-a_{ji}, b_{ij}=-b_{ji}, \forall i, j=1, \ldots, n$, thus proving the matrix of $S$ belongs to $\mathfrak g_{(-1,-1)}$ (see (\ref{eq:(-1-1)-algebra})). $\blacksquare$

\begin{defin}
\label{teor:ae-baselocal}
Let $(M,J,g)$ be a $(J^2=\pm1)$-metric manifold of dimension $2n$ and let $U\subseteq M$be an open subset.  A family $(X_1, \ldots, X_n, Y_1, \ldots, Y_n)$, $X_1, \ldots, X_n, Y_1, \ldots, Y_n \in \mathfrak{X} (U)$ is called an adapted local frame to the $G_{(\alpha,\varepsilon)}$-structure defined by $(J,g)$ in $U$ if $(X_1(p), \ldots, X_n(p), Y_1(p),\ldots Y_n(p)) \in {\mathcal C_{(\alpha,\varepsilon)}}_p, \forall p \in  U$.

Its dual local frame is the family $(\eta_1, \ldots, \eta_n, \omega_1, \ldots \omega_n)$, $\eta_1, \ldots, \eta_n, \omega_1, \ldots, \omega_n \in \bigwedge^1 (U)$, satisfying
\[
\eta_i(X_j)= \omega_{i} (Y_i)=\delta_{ij}, \quad \eta_i (Y_j)=\omega_i(X_j)=0, \quad \forall i, j =1, \ldots, n.
\]
\end{defin}

In the case $\varepsilon=1$ the dual frame is given by
\[
\eta_i (X)= g(X_i,X), \quad \omega_i (X)= g(Y_i,X), \quad i=1, \ldots, n, \quad \forall X\in \mathfrak{X} (U),
\]
while in the case  $\varepsilon=-1$ it is given by
\[
\eta_i (X)= g(Y_i,X), \quad \omega_i (X)= g(X_i,X), \quad i=1, \ldots, n,\quad \forall X\in \mathfrak{X} (U).
\]

\medskip

The following result allow to obtain a local basis of section of the adjoint bundle:

\begin{prop}
\label{teor:base-adjunto-ae-estructura}
Let $(M,J,g)$ be a $(J^2=\pm1)$-metric manifold of dimension $2n$ and let $U\subseteq M$ be an open subset.  Let $(X_1, \ldots, X_n, Y_1, \ldots, Y_n)$ be an adapted local frame to the $G_{(\alpha,\varepsilon)}$-structure defined by $(J,g)$ in $U$ and let $(\eta_1, \ldots, \eta_n, \omega_1, \ldots \omega_n)$ be its dual local frame.

\begin{enumerate}
\renewcommand*{\theenumi}{\roman{enumi})}
\renewcommand*{\labelenumi}{\theenumi}

\item If $(J,g)$ is $(-1,1)$-structure then
\[
\left\{
\begin{array}{c}
S_{ab}=\eta_b \otimes X_a -\eta_a \otimes X_b +\omega_b \otimes Y_a -\omega_a\otimes Y_b\\
S'_{ab}=\eta_b \otimes Y_a +\eta_a \otimes Y_b -\omega_b \otimes X_a -\omega_a\otimes X_b\\
\end{array} \colon a, b \in \{1, \ldots, n\}, 1 \leq a <b \leq n
\right\}
\]
is a local basis of sections in $U$ of the adjoint bundle $\mathrm{ad}  \mathcal C_{(-1,1)} = (\mathcal C_{(-1,1)} \times \mathfrak g_{(-1,1)})/G_{(-1,1)}$.

\item If  $(J,g)$ is a $(-1,-1)$-structure then
\[
\left\{
\begin{array}{c}
S_{ab}=\eta_b \otimes X_a -\eta_a \otimes X_b +\omega_b \otimes Y_a -\omega_a\otimes Y_b\\
S'_{ab}=\eta_b \otimes Y_a -\eta_a \otimes Y_b -\omega_b \otimes X_a +\omega_a\otimes X_b\\
\end{array} \colon a, b \in \{1, \ldots, n\}, 1 \leq a <b \leq n
\right\}
\]
is a local basis of sections in $U$ of the adjoint bundle   $\mathrm{ad}  \mathcal C_{(-1,-1)} = (\mathcal C_{(-1,-1)} \times \mathfrak g_{(-1,-1)})/G_{(-1,-1)}$.

\item If $(J,g)$ is a $(1,1)$-structure then
\[
\left\{
\begin{array}{c}
S_{ab}=\eta_b \otimes X_a -\eta_a \otimes X_b \\
S'_{ab}=\omega_b \otimes Y_a -\omega_a \otimes Y_b \\
\end{array} \colon a, b \in \{1, \ldots, n\}, 1 \leq a <b \leq n
\right\}
\]
is a local basis of sections in $U$ of the adjoint bundle $\mathrm{ad}  \mathcal C_{(1,1)} = (\mathcal C_{(1,1)} \times \mathfrak g_{(1,1)})/G_{(1,1)}$.

\item If $(J,g)$ is $(1,-1)$-structure then
\[
\left\{
S_{ab}=\eta_b \otimes X_a -\omega_a \otimes Y_b \\ \colon a, b \in \{1, \ldots, n\}
\right\}
\]
is a local basis of sections in $U$ of the adjoint bundle
$\mathrm{ad}  \mathcal C_{(1,-1)} = (\mathcal C_{(1,-1)} \times \mathfrak g_{(1,-1)})/G_{(1,-1)}$.
\end{enumerate}
\end{prop}

{\bf Proof.} Trivial, taking into account Proposition \ref{teor:adjunto-ae-estructura}. $\blacksquare$

\bigskip

Then, we can prove the main Theorem of this Section:

\bigskip

\textbf{Proof of Theorem \ref{teor:bienadaptada-ae-estructura}.}
We must prove that the derivation law $\nabla^{\mathrm{w}}$ is characterized by parallelizing $g$ and $J$ and by formula (\ref{eq:welladapted}). Theorem \ref{teor:metodo} shows that the well adapted connection is characterized as the unique natural connection satisfying $ \mathrm{trace}\, (S \circ i_X \circ \mathrm{T}^{\mathrm{w}}) = 0 $, for all section $S$ of the adjoint bundle, where $\mathrm{T}^{\mathrm{w}}$ is the torsion tensor of the derivation law  $\nabla ^{\mathrm{w}}$. Sections of the adjoint bundle have been characterized in  Proposition \ref{teor:adjunto-ae-estructura}. Local basis of sections of the adjoint bundle have been obtained in the previous Proposition \ref{teor:base-adjunto-ae-estructura}. Then, combining all the above results we will able to prove the Theorem.

As in the proof of Proposition \ref{teor:adjunto-ae-estructura}, we will distinguish  four cases. We begin with $\alpha =-1$ and the two subcases $\varepsilon =\pm 1$ and after that we will study the case $\alpha =1$ and the corresponding two subcases. As our proof will be local, we assume that  $(X_1, \ldots, X_n, Y_1, \ldots, Y_n)$ is an adapted local frame to the $G_{(\alpha,\varepsilon)}$-structure defined by $(J,g)$ in an open subset $U$ of $M$.
\medskip

Let $\alpha=-1$ and let us denote $V_1=<X_1, \ldots, X_n>$, $V_2=<Y_1,\ldots, Y_n > = <JX_1, \ldots, JX_n> = JV_1$. For any two vector fields   $X, Z$ in $U$ there exist $X^1, X^2, Z^1, Z^2 \in V_1$ such that
\[
X= X^1+JX^2, \quad Z =Z^1+JZ^2, \quad JX = -X^2 +JX^1, \quad JZ = -Z^2+JZ^1,
\]
thus obtaining
\begin{eqnarray}
g(\mathrm{T}^{\mathrm{w}}(X,Y),Z)-g(\mathrm{T}^{\mathrm{w}}(Z,Y),X)
                                                 &=& g(\mathrm{T}^{\mathrm{w}}(X^1,Y),Z^1)+g(\mathrm{T}^{\mathrm{w}}(JX^2,Y),Z^1) \nonumber\\
                                                 &+&g(\mathrm{T}^{\mathrm{w}}(X^1,Y),JZ^2)+g(\mathrm{T}^{\mathrm{w}}(JX^2,Y),JZ^2) \nonumber\\
                                                 &-& g(\mathrm{T}^{\mathrm{w}}(Z^1,Y),X^1)-g(\mathrm{T}^{\mathrm{w}}(JZ^2,Y),X^1) \nonumber\\
                                                 &-&g(\mathrm{T}^{\mathrm{w}}(Z^1,Y),JX^2)+g(\mathrm{T}^{\mathrm{w}}(JZ^2,Y),JX^2), \label{eq:functorial-11} \\
-\varepsilon(g(\mathrm{T}^{\mathrm{w}}(JX,Y),JZ)-g(\mathrm{T}^{\mathrm{w}}(JZ,Y),JX))
                                                 &=& -\varepsilon \biggl( g(\mathrm{T}^{\mathrm{w}}(X^2,Y),Z^2)-g(\mathrm{T}^{\mathrm{w}}(JX^1,Y),Z^2) \nonumber\\
                                                 &-&g(\mathrm{T}^{\mathrm{w}}(X^2,Y),JZ^1)+g(\mathrm{T}^{\mathrm{w}}(JX^1,Y),JZ^1) \nonumber\\
                                                 &-& g(\mathrm{T}^{\mathrm{w}}(Z^2,Y),X^2)+g(\mathrm{T}^{\mathrm{w}}(JZ^1,Y),X^2) \nonumber\\
                                                 &+&g(\mathrm{T}^{\mathrm{w}}(Z^2,Y),JX^1)-g(\mathrm{T}^{\mathrm{w}}(JZ^1,Y),JX^1) \biggl). \label{eq:functorial-12}
\end{eqnarray}

If $(\alpha,\varepsilon)=(-1,1)$, taking into account  Proposition \ref{teor:base-adjunto-ae-estructura} $i)$, a local basis of sections of $\mathrm{ad} \mathcal C_{(-1,1)}$ is
\[
\left\{
\begin{array}{c}
S_{ab}=\eta_b \otimes X_a -\eta_a \otimes X_b +\omega_b \otimes Y_a -\omega_a\otimes Y_b\\
S'_{ab}=\eta_b \otimes Y_a +\eta_a \otimes Y_b -\omega_b \otimes X_a -\omega_a\otimes X_b\\
\end{array} \colon a, b \in \{1, \ldots, n\}, 1 \leq a <b \leq n
\right\}.
\]
Then, by Theorem \ref{teor:metodo}, $\nabla^{\mathrm{w}}$ is the unique natural derivation law satisfying
\[
\mathrm{trace}\, (S_{ab} \circ i_Y \circ \mathrm{T}^{\mathrm{w}}) = 0, \quad   \mathrm{trace}\, (S'_{ab} \circ i_Y \circ \mathrm{T}^{\mathrm{w}}) = 0, \quad \forall Y\in \mathfrak{X} (M),
\ {\rm for}\ {\rm all}\ a,b\ {\rm with}\  1\leq a <b \leq n.\]

Taking $Y\in \mathfrak{X} (M)$ and $a,b$ with $1\leq a <b \leq n$ one has
\begin{eqnarray*}
\mathrm{trace}\, (S_{ab} \circ i_Y \circ \mathrm{T}^{\mathrm{w}}) &=& \sum_{i=1}^n \eta_i (S_{ab} (\mathrm{T}^{\mathrm{w}}(Y,X_i)))+ \sum_{i=1}^n \omega_i (S_{ab}(\mathrm{T}^{\mathrm{w}}(Y,Y_i)))\\
                                                       &=& \eta_b (\mathrm{T}^{\mathrm{w}}(Y, X_a))-\eta_a(\mathrm{T}^{\mathrm{w}}(Y,X_b))+\omega_b (\mathrm{T}^{\mathrm{w}}(Y,Y_a))-\omega_a(\mathrm{T}^{\mathrm{w}}(Y,Y_b))\\
                                                       &=& g(\mathrm{T}^{\mathrm{w}}(Y, X_a),X_b)-g(\mathrm{T}^{\mathrm{w}}(Y,X_b),X_a)+g (\mathrm{T}^{\mathrm{w}}(Y,Y_a),Y_b)-g(\mathrm{T}^{\mathrm{w}}(Y,Y_b),Y_a),\\
\mathrm{trace}\, (S'_{ab} \circ i_X \circ \mathrm{T}^{\mathrm{w}}) &=& \sum_{i=1}^n \eta_i (S'_{ab} (\mathrm{T}^{\mathrm{w}}(Y,X_i)))+ \sum_{i=1}^n \omega_i (S'_{ab}(\mathrm{T}^{\mathrm{w}}(Y,Y_i)))\\
                                                       &=&- \omega_b (\mathrm{T}^{\mathrm{w}}(Y, X_a))-\omega_a(\mathrm{T}^{\mathrm{w}}(Y,X_b))+\eta_b (\mathrm{T}^{\mathrm{w}}(Y,Y_a))+\eta_a(\mathrm{T}^{\mathrm{w}}(Y,Y_b))\\
                                                       &=& - g (\mathrm{T}^{\mathrm{w}}(Y, X_a), Y_b)-g(\mathrm{T}^{\mathrm{w}}(Y,X_b),Y_a)+g (\mathrm{T}^{\mathrm{w}}(Y,Y_a),X_b)+g(\mathrm{T}^{\mathrm{w}}(Y,Y_b), X_a).
\end{eqnarray*}
From the above conditions one deduces:
\begin{eqnarray*}
g(\mathrm{T}^{\mathrm{w}}(Y, X_a),X_b)-g(\mathrm{T}^{\mathrm{w}}(Y,X_b),X_a)+g (\mathrm{T}^{\mathrm{w}}(Y,Y_a),Y_b)-g(\mathrm{T}^{\mathrm{w}}(Y,Y_b),Y_a) &=&0,\\
- g (\mathrm{T}^{\mathrm{w}}(Y, X_a), Y_b)-g(\mathrm{T}^{\mathrm{w}}(Y,X_b),Y_a)+g (\mathrm{T}^{\mathrm{w}}(Y,Y_a),X_b)+g(\mathrm{T}^{\mathrm{w}}(Y,Y_b), X_a)&=&0,
\end{eqnarray*}
for all $a,b$ with $ 1\leq a < b\leq n$. Given $X^1, Z^1 \in V_1$ one has
\begin{eqnarray}
g(\mathrm{T}^{\mathrm{w}}(X^1, Y),Z^1)-g(\mathrm{T}^{\mathrm{w}}(Z^1,Y),X^1)+g (\mathrm{T}^{\mathrm{w}}(JX^1,Y),JZ^1)-g(\mathrm{T}^{\mathrm{w}}(JZ^1,Y),JX^1) &=&0, \label{eq:functorial-1+11}\\
-g(\mathrm{T}^{\mathrm{w}}(X^1, Y),JZ^1)-g(\mathrm{T}^{\mathrm{w}}(Z^1,Y),JX^1)+g (\mathrm{T}^{\mathrm{w}}(JX^1,Y),Z^1)+g(\mathrm{T}^{\mathrm{w}}(JZ^1,Y),X^1) &=&0. \label{eq:functorial-1+12}
\end{eqnarray}
And given $X^1, X^2, Z^1, Z^2 \in V_1$, from equation  (\ref{eq:functorial-1+11}), one  obtains
\begin{eqnarray*}
g(\mathrm{T}^{\mathrm{w}}(X^1,Y),Z^1)-g(\mathrm{T}^{\mathrm{w}}(Z^1,Y),X^1)&=& -g(\mathrm{T}^{\mathrm{w}}(JX^1,Y),JZ^1)+g(\mathrm{T}^{\mathrm{w}}(JZ^1,Y),JX^1),\\
g(\mathrm{T}^{\mathrm{w}}(JX^2,Y),JZ^2)-g(\mathrm{T}^{\mathrm{w}}(JZ^2,Y),JX^2)&=& -g(\mathrm{T}^{\mathrm{w}}(X^2,Y),Z^2)+g(\mathrm{T}^{\mathrm{w}}(Z^2,Y),X^2),
\end{eqnarray*}
while from equation (\ref{eq:functorial-1+12}) one obtains
\begin{eqnarray*}
g(\mathrm{T}^{\mathrm{w}}(X^1,Y),JZ^2)-g(\mathrm{T}^{\mathrm{w}}(JZ^2,Y),X^1)&=& -g(\mathrm{T}^{\mathrm{w}}(Z^2,Y),JX^1)+g(\mathrm{T}^{\mathrm{w}}(JX^1,Y),Z^2),\\
g(\mathrm{T}^{\mathrm{w}}(JX^2,Y),Z^1)-g(\mathrm{T}^{\mathrm{w}}(Z^1,Y),JX^2)&=& -g(\mathrm{T}^{\mathrm{w}}(X^2,Y),JZ^1)+g(\mathrm{T}^{\mathrm{w}}(JZ^1,Y),X^2),
\end{eqnarray*}
These last equalities combined with (\ref{eq:functorial-11}) and  (\ref{eq:functorial-12}) give the expected relation
\[
g(\mathrm{T}^{\mathrm{w}}(X,Y),Z)-g(\mathrm{T}^{\mathrm{w}}(Z,Y),X)= -(g(\mathrm{T}^{\mathrm{w}}(JX,Y),JZ)-g(\mathrm{T}^{\mathrm{w}}(JZ,Y),JX)), \quad \forall X, Y, Z \in \mathfrak{X} (M).
\]
which is formula (\ref{eq:welladapted}) in the case $\varepsilon=1$.

Observe that last equation in the case $X=X^1$ and $Z=Z^1$ reads as
\[
g(\mathrm{T}^{\mathrm{w}}(X^1,Y),Z^1)-g(\mathrm{T}^{\mathrm{w}}(Z^1,Y),X^1)+g(\mathrm{T}^{\mathrm{w}}(JX^1,Y),JZ^1)-g(\mathrm{T}^{\mathrm{w}}(JZ^1,Y),JX^1)=0,
\]
i.e, coincides with formula (\ref{eq:functorial-1+11}), while in the case $X=-X^1$ and $Z=JZ^1$ reads as
\[
-g(\mathrm{T}^{\mathrm{w}}(X^1,Y),JZ^1)+g(\mathrm{T}^{\mathrm{w}}(JZ^1,Y),X^1)+g(\mathrm{T}^{\mathrm{w}}(JX^1,Y),Z^1)-g(\mathrm{T}^{\mathrm{w}}(Z^1,Y),JX^1)=0
\]
and thus coincides with formula (\ref{eq:functorial-1+12}).
\medskip

If $(\alpha,\varepsilon)=(-1,-1)$ taking into account Proposition \ref{teor:base-adjunto-ae-estructura} $ii)$, a local basis of sections of $\mathrm{ad} \mathcal C_{(-1,-1)}$ is
\[
\left\{
\begin{array}{c}
S_{ab}=\eta_b \otimes X_a -\eta_a \otimes X_b +\omega_b \otimes Y_a -\omega_a\otimes Y_b\\
S'_{ab}=\eta_b \otimes Y_a -\eta_a \otimes Y_b -\omega_b \otimes X_a +\omega_a\otimes X_b\\
\end{array} \colon a, b \in \{1, \ldots, n\}, 1 \leq a <b \leq n
\right\}
\]
Then, by Theorem \ref{teor:metodo} $\nabla^{\mathrm{w}}$ is the unique natural derivation law satisfying
\[
\mathrm{trace}\, (S_{ab} \circ i_Y \circ \mathrm{T}^{\mathrm{w}}) = 0, \quad   \mathrm{trace}\, (S'_{ab} \circ i_Y \circ \mathrm{T}^{\mathrm{w}}) = 0, \quad \forall Y\in \mathfrak{X} (M),
\ {\rm for}\ {\rm all}\ a,b\ {\rm with}\  1\leq a <b \leq n.\]

Taking $Y\in \mathfrak{X} (M)$ and $a,b$ with $1\leq a <b \leq n$ one has
\begin{eqnarray*}
\mathrm{trace}\, (S_{ab} \circ i_Y \circ \mathrm{T}^{\mathrm{w}}) &=& \sum_{i=1}^n \eta_i (S_{ab} (\mathrm{T}^{\mathrm{w}}(Y,X_i)))+ \sum_{i=1}^n \omega_i (S_{ab}(\mathrm{T}^{\mathrm{w}}(Y,Y_i)))\\
                                                       &=& \eta_b (\mathrm{T}^{\mathrm{w}}(Y, X_a))-\eta_a(\mathrm{T}^{\mathrm{w}}(Y,X_b))+\omega_b (\mathrm{T}^{\mathrm{w}}(Y,Y_a))-\omega_a(\mathrm{T}^{\mathrm{w}}(Y,Y_b))\\
                                                       &=& g(\mathrm{T}^{\mathrm{w}}(Y, X_a),Y_b)-g(\mathrm{T}^{\mathrm{w}}(Y,X_b),Y_a)+g (\mathrm{T}^{\mathrm{w}}(Y,Y_a),X_b)-g(\mathrm{T}^{\mathrm{w}}(Y,Y_b),X_a),
\\
\mathrm{trace}\, (S'_{ab} \circ i_X \circ \mathrm{T}^{\mathrm{w}}) &=& \sum_{i=1}^n \eta_i (S'_{ab} (\mathrm{T}^{\mathrm{w}}(Y,X_i)))+ \sum_{i=1}^n \omega_i (S'_{ab}(\mathrm{T}^{\mathrm{w}}(Y,Y_i)))\\
                                                       &=&- \omega_b (\mathrm{T}^{\mathrm{w}}(Y, X_a))+\omega_a(\mathrm{T}^{\mathrm{w}}(Y,X_b))+\eta_b (\mathrm{T}^{\mathrm{w}}(Y,Y_a))-\eta_a(\mathrm{T}^{\mathrm{w}}(Y,Y_b))\\
                                                       &=& - g (\mathrm{T}^{\mathrm{w}}(Y, X_a), X_b)+g(\mathrm{T}^{\mathrm{w}}(Y,X_b),X_a)+g (\mathrm{T}^{\mathrm{w}}(Y,Y_a),Y_b)-g(\mathrm{T}^{\mathrm{w}}(Y,Y_b), Y_a).
\end{eqnarray*}
From the above conditions one deduces:
\begin{eqnarray*}
g(\mathrm{T}^{\mathrm{w}}(Y, X_a),Y_b)-g(\mathrm{T}^{\mathrm{w}}(Y,X_b),Y_a)+g (\mathrm{T}^{\mathrm{w}}(Y,Y_a),X_b)-g(\mathrm{T}^{\mathrm{w}}(Y,Y_b),X_a) &=&0,\\
g (\mathrm{T}^{\mathrm{w}}(Y, X_a), X_b)+g(\mathrm{T}^{\mathrm{w}}(Y,X_b),X_a)+g (\mathrm{T}^{\mathrm{w}}(Y,Y_a),Y_b)-g(\mathrm{T}^{\mathrm{w}}(Y,Y_b), Y_a)&=&0,
\end{eqnarray*}
for all $a,b$ with  $1\leq a < b\leq n$. Given $X^1, Z^1 \in V_1$ one has
\begin{eqnarray}
g(\mathrm{T}^{\mathrm{w}}(X^1, Y),JZ^1)-g(\mathrm{T}^{\mathrm{w}}(Z^1,Y),JX^1)+g (\mathrm{T}^{\mathrm{w}}(JX^1,Y),Z^1)-g(\mathrm{T}^{\mathrm{w}}(JZ^1,Y),X^1) &=&0, \label{eq:functorial-1+13}\\
-g(\mathrm{T}^{\mathrm{w}}(X^1, Y),Z^1)+g(\mathrm{T}^{\mathrm{w}}(Z^1,Y),X^1)+g (\mathrm{T}^{\mathrm{w}}(JX^1,Y),JZ^1)-g(\mathrm{T}^{\mathrm{w}}(JZ^1,Y),X^1) &=&0. \label{eq:functorial-1+14}
\end{eqnarray}
And given $X^1, X^2, Z^1, Z^2 \in V_1$, from equation (\ref{eq:functorial-1+13}), one obtains
\begin{eqnarray*}
g(\mathrm{T}^{\mathrm{w}}(X^1,Y),JZ^2)-g(\mathrm{T}^{\mathrm{w}}(JZ^2,Y),X^1)&=& g(\mathrm{T}^{\mathrm{w}}(Z^2,Y),JX^1)-g(\mathrm{T}^{\mathrm{w}}(JX^1,Y),Z^2),\\
g(\mathrm{T}^{\mathrm{w}}(JX^2,Y),Z^1)-g(\mathrm{T}^{\mathrm{w}}(Z^1,Y),JX^2)&=& -g(\mathrm{T}^{\mathrm{w}}(X^2,Y),JZ^1)+g(\mathrm{T}^{\mathrm{w}}(JZ^1,Y),X^2),
\end{eqnarray*}
while from equation (\ref{eq:functorial-1+14}) one obtains
\begin{eqnarray*}
g(\mathrm{T}^{\mathrm{w}}(X^1,Y),Z^1)-g(\mathrm{T}^{\mathrm{w}}(Z^1,Y),X^1)&=& g(\mathrm{T}^{\mathrm{w}}(JX^1,Y),JZ^1)-g(\mathrm{T}^{\mathrm{w}}(JZ^1,Y),JX^1),\\
g(\mathrm{T}^{\mathrm{w}}(JX^2,Y),JZ^2)-g(\mathrm{T}^{\mathrm{w}}(JZ^2,Y),JX^2)&=& g(\mathrm{T}^{\mathrm{w}}(Z^2,Y),Z^2)-g(\mathrm{T}^{\mathrm{w}}(Z^2,Y),X^2),
\end{eqnarray*}
These last equalities combined with (\ref{eq:functorial-11}) and (\ref{eq:functorial-12}) give the expected relation
\[
g(\mathrm{T}^{\mathrm{w}}(X,Y),Z)-g(\mathrm{T}^{\mathrm{w}}(Z,Y),X)= g(\mathrm{T}^{\mathrm{w}}(JX,Y),JZ)-g(\mathrm{T}^{\mathrm{w}}(JZ,Y),JX), \quad \forall X, Y, Z \in \mathfrak{X} (M).
\]
which is formula (\ref{eq:welladapted}) in the case $\varepsilon=-1$.

Observe that the last equation in the case $X=X^1$ and $Z=JZ^1$ reads as
\[
g(\mathrm{T}^{\mathrm{w}}(X^1,Y),JZ^1)-g(\mathrm{T}^{\mathrm{w}}(JZ^1,Y),X^1)+g(\mathrm{T}^{\mathrm{w}}(JX^1,Y),Z^1)-g(\mathrm{T}^{\mathrm{w}}(Z^1,Y),JX^1)=0,
\]
i.e., coincides with formula (\ref{eq:functorial-1+13}), while in the case $X=-X^1$ and  $Z=Z^1$ reads as
\[
-g(\mathrm{T}^{\mathrm{w}}(X^1,Y),Z^1)+g(\mathrm{T}^{\mathrm{w}}(Z^1,Y),X^1)+g(\mathrm{T}^{\mathrm{w}}(JX^1,Y),JZ^1)-g(\mathrm{T}^{\mathrm{w}}(JZ^1,Y),JX^1)=0
\]
and thus coincides with formula (\ref{eq:functorial-1+14}).
\medskip

Let $\alpha=1$. Let us denote $T^+_J(U)=<X_1, \ldots, X_n>$, $T^-_J(U)=<Y_1,\ldots, Y_n >=<JX_1, \ldots, JX_n>$, and
\[
J^+=\frac{1}{2} (Id+J), \quad J^-=\frac{1}{2}(Id-J).
\]
Given two vector fields $X, Z$ in $U$ one has
\[
X=J^+X + J^-X, \quad Z=J^+Z+J^-Z, \quad JX = J^+X-J^-X, \quad JZ =J^+Z-J^-Z,
\]
thus obtaining
\begin{eqnarray}
g(\mathrm{T}^{\mathrm{w}}(X,Y),Z)-g(\mathrm{T}^{\mathrm{w}}(Z,Y),X)&=& g(\mathrm{T}^{\mathrm{w}}( J^+X + J^-X,Y),J^+Z + J^-Z)\nonumber \\
                                                 &-&g(\mathrm{T}^{\mathrm{w}}(J^+Z + J^-Z,Y),J^+X + J^-X) \nonumber\\
                                                 &=& g(\mathrm{T}^{\mathrm{w}}(J^+X,Y),J^+Z)+g(\mathrm{T}^{\mathrm{w}}(J^-X,Y),J^+Z) \nonumber\\
                                                 &+&g(\mathrm{T}^{\mathrm{w}}(J^+X,Y),J^-Z)+g(\mathrm{T}^{\mathrm{w}}(J^-X,Y),J^-Z) \nonumber\\
                                                 &-& g(\mathrm{T}^{\mathrm{w}}(J^+Z,Y),J^+X)-g(\mathrm{T}^{\mathrm{w}}(J^-Z,Y),J^+X) \nonumber\\
                                                 &-&g(\mathrm{T}^{\mathrm{w}}(J^+Z,Y),J^-X)-g(\mathrm{T}^{\mathrm{w}}(J^-Z,Y),J^-X), \label{eq:functorial+11}
\end{eqnarray}
\begin{eqnarray}             
                                                 -\varepsilon(g(\mathrm{T}^{\mathrm{w}}(JX,Y),JZ)-g(\mathrm{T}^{\mathrm{w}}(JZ,Y),JX))&=& -\varepsilon\biggl(g(\mathrm{T}^{\mathrm{w}}( J^+X - J^-X,Y),J^+Z - J^-Z)\nonumber  \\
                                                                  &-& g(\mathrm{T}^{\mathrm{w}}(J^+Z - J^-Z,Y),J^+X - J^-X) \biggl)\nonumber\\
                                                 &=& -\varepsilon \biggl( g(\mathrm{T}^{\mathrm{w}}(J^+X,Y),J^+Z)-g(\mathrm{T}^{\mathrm{w}}(J^-X,Y),J^+Z) \nonumber\\
                                                 &-&g(\mathrm{T}^{\mathrm{w}}(J^+X,Y),J^-Z)+g(\mathrm{T}^{\mathrm{w}}(J^-X,Y),J^-Z) \nonumber\\
                                                 &-& g(\mathrm{T}^{\mathrm{w}}(J^+Z,Y),J^+X)+g(\mathrm{T}^{\mathrm{w}}(J^-Z,Y),J^+X) \nonumber\\
                                                 &+&g(\mathrm{T}^{\mathrm{w}}(J^+Z,Y),J^-X)-g(\mathrm{T}^{\mathrm{w}}(J^-Z,Y),J^-X) \biggl). \label{eq:functorial+12}
\end{eqnarray}

If $(\alpha,\varepsilon)=(1,1)$ taking into account Proposition \ref{teor:base-adjunto-ae-estructura} $iii)$, a local basis of sections of $\mathrm{ad} \mathcal C_{(1,1)}$ is
\[
\left\{
\begin{array}{c}
S_{ab}=\eta_b \otimes X_a -\eta_a \otimes X_b \\
S'_{ab}=\omega_b \otimes Y_a -\omega_a \otimes Y_b
\end{array} \colon a, b \in \{1, \ldots, n\}, 1 \leq a <b \leq n
\right\}
\]
Then, by Theorem \ref{teor:metodo}  $\nabla^{\mathrm{w}}$ is the unique natural derivation law satisfying
\[
\mathrm{trace}\, (S_{ab} \circ i_Y \circ \mathrm{T}^{\mathrm{w}}) = 0, \quad   \mathrm{trace}\, (S'_{ab} \circ i_Y \circ \mathrm{T}^{\mathrm{w}}) = 0, \quad \forall Y\in \mathfrak{X} (M),
\ {\rm for}\ {\rm all}\ a,b\ {\rm with}\  1\leq a <b \leq n.\]

Taking $Y\in \mathfrak{X} (M)$ and $a,b$ with $1\leq a <b \leq n$ one has
\begin{eqnarray*}
\mathrm{trace}\, (S_{ab} \circ i_Y \circ \mathrm{T}^{\mathrm{w}}) &=& \sum_{i=1}^n \eta_i (S_{ab} (\mathrm{T}^{\mathrm{w}}(Y,X_i)))+ \sum_{i=1}^n \omega_i (S_{ab}(\mathrm{T}^{\mathrm{w}}(Y,Y_i)))\\
                                                       &=& \eta_b(\mathrm{T}^{\mathrm{w}}(Y,X_a))-\eta_a(\mathrm{T}^{\mathrm{w}}(Y,X_b))= g(\mathrm{T}^{\mathrm{w}}(Y,X_a),X_b)-g(\mathrm{T}^{\mathrm{w}}(Y,X_b),X_a),
\\
\mathrm{trace}\, (S'_{ab} \circ i_X \circ \mathrm{T}^{\mathrm{w}}) &=& \sum_{i=1}^n \eta_i (S'_{ab} (\mathrm{T}^{\mathrm{w}}(Y,X_i)))+ \sum_{i=1}^n \omega_i (S'_{ab}(\mathrm{T}^{\mathrm{w}}(Y,Y_i)))\\
                                                       &=& \omega_b(\mathrm{T}^{\mathrm{w}}(Y,Y_a))-\omega_a(\mathrm{T}^{\mathrm{w}}(Y,Y_b))= g(\mathrm{T}^{\mathrm{w}}(Y,Y_a),Y_b)-g (\mathrm{T}^{\mathrm{w}}(Y,Y_b),Y_a).
\end{eqnarray*}
From the above conditions one deduces
\[
g(\mathrm{T}^{\mathrm{w}}(Y,X_a),X_b)-g(\mathrm{T}^{\mathrm{w}}(Y,X_b),X_a)) =0, \quad
g(\mathrm{T}^{\mathrm{w}}(Y,Y_a),Y_b)-g (\mathrm{T}^{\mathrm{w}}(Y,Y_b),Y_a)=0,
\]
for all $a,b$ with $1\leq a < b\leq n$. Given vector fields $X, Z$ in $U$ one has
\begin{eqnarray}
g(\mathrm{T}^{\mathrm{w}}(J^+X,Y),J^+Z)-g(\mathrm{T}^{\mathrm{w}}(J^+Z,Y),J^+X) &=&0, \label{eq:functorial+1+11}\\
g(\mathrm{T}^{\mathrm{w}}(J^-X,Y),J^-Z)-g(\mathrm{T}^{\mathrm{w}}(J^-Z,Y),J^-X)&=&0. \label{eq:functorial+1+12}
\end{eqnarray}
By applying  (\ref{eq:functorial+1+11}) and (\ref{eq:functorial+1+12}) to (\ref{eq:functorial+11}) one deduces
\begin{eqnarray*}
g(\mathrm{T}^{\mathrm{w}}(X,Y),Z)-g(\mathrm{T}^{\mathrm{w}}(Z,Y),X) &=& g(\mathrm{T}^{\mathrm{w}}(J^-X,Y),J^+Z)+g(\mathrm{T}^{\mathrm{w}}(J^+X,Y),J^-Z)\\
                                                 &-&g(\mathrm{T}^{\mathrm{w}}(J^-Z,Y),J^+X) -g(\mathrm{T}^{\mathrm{w}}(J^+Z,Y),J^-X),
\end{eqnarray*}
while applying (\ref{eq:functorial+1+11}) and (\ref{eq:functorial+1+12}) to (\ref{eq:functorial+12}) one obtains
\begin{eqnarray*}
-(g(\mathrm{T}^{\mathrm{w}}(JX,Y),JZ)-g(\mathrm{T}^{\mathrm{w}}(JZ,Y),JX))&=& g(\mathrm{T}^{\mathrm{w}}(J^-X,Y),J^+Z)+g(\mathrm{T}^{\mathrm{w}}(J^+X,Y),J^-Z)\\
                                                 &-&g(\mathrm{T}^{\mathrm{w}}(J^-Z,Y),J^+X) -g(\mathrm{T}^{\mathrm{w}}(J^+Z,Y),J^-X),
\end{eqnarray*}
thus proving
\[
g(\mathrm{T}^{\mathrm{w}}(X,Y),Z)-g(\mathrm{T}^{\mathrm{w}}(Z,Y),X)= -(g(\mathrm{T}^{\mathrm{w}}(JX,Y),JZ)-g(\mathrm{T}^{\mathrm{w}}(JZ,Y),JX)), \quad \forall X, Y, Z \in \mathfrak{X} (M).
\]

Observe that in the above equation when  $(X,Z)$ are changed by $(J^+X,J^+Z)$ one obtains
\begin{eqnarray*}
g(\mathrm{T}^{\mathrm{w}}(J^+X,Y),J^+Z)-g(\mathrm{T}^{\mathrm{w}}(J^+Z,Y),J^+X)&=&-g(\mathrm{T}^{\mathrm{w}}(J^+X,Y),J^+Z)+g(\mathrm{T}^{\mathrm{w}}(J^+Z,Y),J^+X),\\
g(\mathrm{T}^{\mathrm{w}}(J^+X,Y),J^+Z) &=& g(\mathrm{T}^{\mathrm{w}}(J^+Z,Y),J^+X),
\end{eqnarray*}
i.e., one obtains formula (\ref{eq:functorial+1+11}), while changing $(X,Z)$ by $(J^-X,J^-Z)$ one has
\begin{eqnarray*}
g(\mathrm{T}^{\mathrm{w}}(J^-X,Y),J^-Z)-g(\mathrm{T}^{\mathrm{w}}(J^-Z,Y),J^-X)&=& - g(\mathrm{T}^{\mathrm{w}}(J^-X,Y),J^-Z)+g(\mathrm{T}^{\mathrm{w}}(J^-Z,Y),J^-X),\\
g(\mathrm{T}^{\mathrm{w}}(J^-X,Y),J^-Z)&=& g(\mathrm{T}^{\mathrm{w}}(J^-Z,Y),J^-X),
\end{eqnarray*}
which is formula (\ref{eq:functorial+1+12}).
\medskip

If $(\alpha,\varepsilon)=(1,-1)$ taking into account Proposition \ref{teor:base-adjunto-ae-estructura} $iv)$, a local basis of sections of $\mathrm{ad} \mathcal C_{(1,-1)}$ is
\[
\left\{
S_{ab}=\eta_b \otimes X_a -\omega_a \otimes Y_b \\ \colon a, b \in \{1, \ldots, n\}
\right\}.
\]
Then, by Theorem \ref{teor:metodo} $\nabla^{\mathrm{w}}$is the unique natural derivation law satisfying
\[
\mathrm{trace}\, (S_{ab} \circ i_Y \circ \mathrm{T}^{\mathrm{w}}) = 0, \quad \forall Y\in \mathfrak{X} (M),
\ {\rm for}\ {\rm all}\ a,b\ {\rm with}\  1\leq a <b \leq n.\]

Taking $Y\in \mathfrak{X} (M)$ and $a,b$ with $1\leq a,b \leq n$ one has
\begin{eqnarray*}
\mathrm{trace}\, (S_{ab} \circ i_Y \circ \mathrm{T}^{\mathrm{w}}) &=& \sum_{i=1}^n \eta_i (S_{ab} (\mathrm{T}^{\mathrm{w}}(Y,X_i)))+ \sum_{i=1}^n \omega_i (S_{ab}(\mathrm{T}^{\mathrm{w}}(Y,Y_i)))\\
                                                       &=& \eta_b(\mathrm{T}^{\mathrm{w}}(Y,X_a))-\omega_a(\mathrm{T}^{\mathrm{w}}(Y,Y_b))=g(\mathrm{T}^{\mathrm{w}}(Y,X_a),Y_b)-g(\mathrm{T}^{\mathrm{w}}(Y,Y_b),X_a),
\end{eqnarray*}
and then
\[
g(\mathrm{T}^{\mathrm{w}}(Y,X_a),Y_b)-g(\mathrm{T}^{\mathrm{w}}(Y,Y_b),X_a)=0, \quad \forall a, b=1, \ldots, n.
\]
Consequently, given two vector fields $X, Z$ in $U$ one has
\begin{equation}
g(\mathrm{T}^{\mathrm{w}}(J^+X,Y),J^-Z)-g(\mathrm{T}^{\mathrm{w}}(J^-Z,Y),J^+X) =0,
\label{eq:functorial+1+13}\\
\end{equation}
By applying (\ref{eq:functorial+1+13}) to (\ref{eq:functorial+11}) one deduces
\begin{eqnarray*}
g(\mathrm{T}^{\mathrm{w}}(X,Y),Z)-g(\mathrm{T}^{\mathrm{w}}(Z,Y),X) &=& g(\mathrm{T}^{\mathrm{w}}(J^+X,Y),J^+Z)+g(\mathrm{T}^{\mathrm{w}}(J^-X,Y),J^-Z)\\
                                                 &-&g(\mathrm{T}^{\mathrm{w}}(J^+Z,Y),J^+X) -g(\mathrm{T}^{\mathrm{w}}(J^+Z,Y),J^+X),
\end{eqnarray*}
while applying (\ref{eq:functorial+1+13}) to (\ref{eq:functorial+12}) one obtains
\begin{eqnarray*}
g(\mathrm{T}^{\mathrm{w}}(JX,Y),JZ)-g(\mathrm{T}^{\mathrm{w}}(JZ,Y),JX)&=& g(\mathrm{T}^{\mathrm{w}}(J^-X,Y),J^-Z)+g(\mathrm{T}^{\mathrm{w}}(J^-X,Y),J^-Z)\\
                                                 &-&g(\mathrm{T}^{\mathrm{w}}(J^+Z,Y),J^+X) -g(\mathrm{T}^{\mathrm{w}}(J^-Z,Y),J^-X),
\end{eqnarray*}
thus proving
\[
g(\mathrm{T}^{\mathrm{w}}(X,Y),Z)-g(\mathrm{T}^{\mathrm{w}}(Z,Y),X)= g(\mathrm{T}^{\mathrm{w}}(JX,Y),JZ)-g(\mathrm{T}^{\mathrm{w}}(JZ,Y),JX), \quad \forall X, Y, Z \in \mathfrak{X} (M).
\]

Observe that in the above equation when  $(X,Z)$ are changed by $(J^+X,J^-Z)$ one obtains
\begin{eqnarray*}
g(\mathrm{T}^{\mathrm{w}}(J^+X,Y),J^-Z)-g(\mathrm{T}^{\mathrm{w}}(J^-Z,Y),J^+X)&=&-g(\mathrm{T}^{\mathrm{w}}(J^+X,Y),J^-Z)+g(\mathrm{T}^{\mathrm{w}}(J^-Z,Y),J^+X),\\
g(\mathrm{T}^{\mathrm{w}}(J^+X,Y),J^-Z) &=& g(\mathrm{T}^{\mathrm{w}}(J^-Z,Y),J^+X),
\end{eqnarray*}
i.e., one obtains formula (\ref{eq:functorial+1+13}). $\blacksquare$

\section{Particularizing the well adapted connection}
\label{sec:particularizingthewelladaptedconnection}

 The expression (\ref{eq:welladapted}) in Theorem \ref{teor:bienadaptada-ae-estructura} is common for the four well adapted connections corresponding to the four classes of $(J^2=\pm1)$-metrics manifolds. We will study them carefully in order to recover connections firstly introduced in the Literature under other denominations.

 \bigskip

\textbf{ Almost para-Hermitian manifolds or $(1,-1)$-structures.}

If $(J,g)$ is an $(1,-1)$-structure,  then the expression (\ref{eq:functorial+1+13})
\[
g(\mathrm{T}^{\mathrm{w}}(J^+X,Y),J^-Z)=g(\mathrm{T}^{\mathrm{w}}(J^-Z,Y),J^+X), \quad \forall X, Y, Z \in \mathfrak{X} (M),
\]
characterize $\nabla^{\mathrm{w}}$. This equation and the conditions $\nabla^{\mathrm{w}} J=0$ and $\nabla^{\mathrm{w}} g=0$ are the characterization of the well adapted connection that we have previously obtained  in \cite[Theor.\ 3.8]{brassov}.

\bigskip

\textbf{ Almost Hermitian manifolds or $(-1,1)$-structures.}

In this case, $\varepsilon =1$,   formula (\ref{eq:welladapted}) reads as
\begin{equation}
g(\mathrm{T}^{\mathrm{w}}(X,Y),Z)-g(\mathrm{T}^{\mathrm{w}}(Z,Y),X) +g(\mathrm{T}^{\mathrm{w}}(JX,Y),JZ)-g(\mathrm{T}^{\mathrm{w}}(JZ,Y),JX)=0, \quad \forall X,Y,Z\in \mathfrak{X} (M).
\label{eq:welladapted-1}
\end{equation}

This property and the parallelization of $J$ and $g$ are used in \cite[Secc.\ 3.1]{valdes} to characterize $\nabla^{\mathrm{w}}$.
By using functorial connections, in \cite[Theor.\ 3.1]{munoz}, the authors had obtained the expression of $\nabla^{\mathrm{w}}$ in terms of complex frames.

\bigskip

In the above two cases, explicit expressions of $\nabla^{\mathrm{w}}$ had been obtained. The situation is different in the other two cases: the expression had been obtained under the name of canonical connections, because the starting point was not that of the well adapted connections. Authors of the ``Bulgarian School" had defined adapted connections which coincide with the well adapted one. Let us see it.

\bigskip

\textbf{ Almost Norden manifolds or $(-1,-1)$-structures and almost product Riemannian manifolds or $(1,1)$-structures.}

If $\varepsilon =-1$ formula (\ref{eq:welladapted}) reads
\begin{equation}
g(\mathrm{T}^{\mathrm{w}}(X,Y),Z)-g(\mathrm{T}^{\mathrm{w}}(Z,Y),X) -g(\mathrm{T}^{\mathrm{w}}(JX,Y),JZ)+g(\mathrm{T}^{\mathrm{w}}(JZ,Y),JX)=0, \quad \forall X,Y,Z\in \mathfrak{X} (M).
\label{eq:welladapted+1}
\end{equation}

Equations (\ref{eq:welladapted-1}) and (\ref{eq:welladapted+1}) are also the conditions defining the well adapted connections in the almost product Riemannian and almost Norden  manifolds, because formula (\ref{eq:welladapted}) depends on $\varepsilon$ but does not depend on $\alpha$. In the almost Norden case, coincides with equations obtained in \cite[Theor. 5]{ganchev-mihova}, and in the almost semi-Riemannian product of signature $(n,n)$ with those of  \cite[Theor.\ 4]{mihova}. In these papers these connections are called canonical connections.

Besides, these cases correspond with  $(J^2=\pm1)$-metric manifolds having $\alpha \varepsilon=1$. The authors of the quoted papers have introduced the canonical connection by using the tensor of type $(1,2)$ defined as the difference of the Levi Civita connections of the metrics $g$ and $\widetilde g$,
where $\widetilde g$ is the twin metric defined as
\[
\widetilde g(X,Y)= g(JX,Y), \quad \forall X, Y \in \mathfrak{X} (M).
\]
The twin metric is semi-Riemannian of signature $(n,n)$ and has a r\^{o}le similar to that of the K\"{a}hler form in the case $\alpha \varepsilon=-1$. 
\bigskip

We finish this section studying the the difference $\nabla^{\mathrm{w}} - \nabla^{\mathrm{g}}$.

\begin{prop}
Let $(M,J,g)$ be a $(J^2=\pm1)$-metric manifold. Then $\nabla^{\mathrm{w}} = \nabla^{\mathrm{g}}$ if and only if
\[
g(\mathrm{T}^{\mathrm{w}}(JX,Y),JZ)=g(\mathrm{T}^{\mathrm{w}}(JZ,Y),JX), \quad \forall X, Y, Z \in \mathfrak{X} (M).
\]

In that case, $(M,J,g)$ is a K\"{a}hler type manifold.
\end{prop}

\textbf{Proof.} The Levi Civita connection $\nabla^{\mathrm{g}}$ is the unique torsionless connection parallelizing $g$. It is easy to prove that these conditions are equivalent to $\nabla^{\mathrm{g}}g=0$ and

\begin{equation*}
g(\mathrm{T}^{\mathrm{g}}(X,Y),Z) = g(\mathrm{T}^{\mathrm{g}}(Z,Y),X), \quad \forall X,Y, Z \in \mathfrak{X} (M).
\end{equation*}

As $\nabla^{\mathrm{w}}$ is an adapted connection one has $\nabla^{\mathrm{w}}g=0$. If $\nabla^{\mathrm{w}}$ satisfies the condition

\begin{equation}
g(\mathrm{T}^{\mathrm{w}}(X,Y),Z) = g(\mathrm{T}^{\mathrm{w}}(Z,Y),X), \quad \forall X,Y, Z \in \mathfrak{X} (M),
\label{eq:ba-ae-torsion2}
\end{equation}

\noindent then $\nabla^{\mathrm{w}}=\nabla^{\mathrm{g}}$. Taking into account formula (\ref{eq:welladapted}), the equality (\ref{eq:ba-ae-torsion2}) is satisfied if and only if

\[
g(\mathrm{T}^{\mathrm{w}}(JX,Y),JZ)=g(\mathrm{T}^{\mathrm{w}}(JZ,Y),JX), \quad \forall X, Y, Z \in \mathfrak{X} (M),
\]

\noindent as we wanted.

Finally, as $\nabla^{\mathrm{w}}$ is an adapted connection one has $\nabla^{\mathrm{w}}J=0$, and then if both connections coincide, the manifold is of K\"{a}hler type because $\nabla^{\mathrm{g}}J=0$, thus finishing the proof. $\blacksquare$

\begin{teor}
\label{Levi Civita2}
 Let $(M,J,g)$ be a $(J^2=\pm1)$-metric manifold. Then $\nabla^{\mathrm{w}} = \nabla^{\mathrm{g}}$ if and only if $(M,J,g)$ is a K\"{a}hler type manifold.
\end{teor}

\textbf{Proof.} The manifold $(M,J,g)$ is a K\"{a}hler type manifold if and only if $\nabla^{\mathrm{g}} J=0$. As the Levi Civita connection of $g$ satisfies $\nabla^{\mathrm{g}} g = 0$ and $\mathrm{T}^{\mathrm{g}}=0$, then the derivation law $\nabla^{\mathrm{g}}$ also satisfies all conditions in Theorem \ref{teor:bienadaptada-ae-estructura}, thus proving the Levi Civita  and the well adapted connections coincide.

Note that the other implication has been proved in the previous proposition, thus finishing the proof. $\blacksquare$

\bigskip

The above results show that the well adapted connection is the most natural extension of the Levi Civita connection to  $(J^2=\pm1)$-metric manifolds.

\section{The Chern connection  of an $(\alpha ,\varepsilon)$-manifold satisfying $\alpha \varepsilon =-1$}
\label{sec:chernconnection}

As we have pointed out in Section~\ref{sec:introduction}, there have been published several papers where the authors were looking for a \textit{canonical} connection in some of the four geometries, generalizing the Levi Civita connection. In the case of $(\alpha,\varepsilon)$-structures with $\alpha\varepsilon=1$ the connections obtained in \cite{ganchev-mihova} and \cite{mihova} coincide with the well adapted connection, as we have seen  in the above Section. In the case $\alpha\varepsilon=-1$ one can define Chern-type connections, which in general do not coincide with the well adapted connection. In this Section we  are going to define Chern-type connections for  $\alpha\varepsilon=-1$, proving that one can not define them in the case  $\alpha\varepsilon=1$, and finally we will characterize when such a connection coincides with the well adapted connection.
\bigskip

The Chern connection was introduced in  \cite{chern}, assuming one has a Hermitian manifold. It also runs in the non-integrable case, because one has:

\begin{teor}[{\cite[Theor. 6.1]{gray}}]
Let  $(M,J,g)$ be an almost Hermitian manifold, i.e., a manifold endowed with a $(-1,1)$-structure.  Then there exists a unique linear connection $\nabla^{\mathrm{c}}$ in  $M$  satisfaying $\nabla^{\mathrm{c}} J=0$, $\nabla^{\mathrm{c}} g=0$ and
\begin{equation*}
\mathrm{T}^{\mathrm{c}}(JX,JY)= -\mathrm{T}^{\mathrm{c}}(X,Y), \quad \forall X, Y \in \mathfrak{X} (M).
\end{equation*}
where $\mathrm{T}^{\mathrm{c}}$ denotes the torsion tensor of $\nabla^{\mathrm{c}}$.
\end{teor}

In the almost para-Hermitian case the connection also runs. In order to prove it, we need the following result:

\begin{lema}\label{teor:chernparahermitian}
Let $(M,J,g)$ be an almost para-Hermitian manifold, i.e, a manifold endowed with a $(1,-1)$-structure.   Let $\nabla$ be any linear connection in $M$ with torsion tensor $\mathrm{T}$. Then, the following conditions are equivalent:
\begin{enumerate}
\renewcommand*{\theenumi}{\roman{enumi})}
\renewcommand*{\labelenumi}{\theenumi}

\item $\mathrm{T}(JX,JY)= \mathrm{T}(X,Y)$, for all vector fields $X, Y$ in $M$.

\item $ \mathrm{T}(X,Y) =0$, for all vector fields  $X \in T^+_J (M)$ and  $Y \in T^-_J(M)$.

\end{enumerate}
\end{lema}

{\bf Proof.}

$ i) \Rightarrow ii)$ Given  $X \in T^+_J (M), Y \in T^-_J(M)$ one has  $JX=X, JY=-Y$.  Then,
\[
\mathrm{T}(JX,JY) = \mathrm{T}(X,Y)
\Rightarrow
-\mathrm{T}(X,Y)=\mathrm{T}(X,Y) \Rightarrow \mathrm{T}(X,Y) =0.
\]

$ii) \Rightarrow i)$  Given two vector fields $X, Y$ in $M$ one has the decompositions
\[
X=J^+X + J^-X, \quad Y = J^+Y-J^-Y, \quad  JX = J^+X -J^-X, \quad JY=J^+Y-J^-Y,
\]
and then taking into account $ii)$ one has
\begin{eqnarray*}
\mathrm{T}(X,Y)&=& \mathrm{T} (J^+X,J^+Y)+\mathrm{T}(J^+X,J^-Y)+\mathrm{T}(J^-X,J^+Y)+\mathrm{T}(J^-X,J^-Y)\\
                                                 &=& \mathrm{T}(J^+X,J^+Y)+\mathrm{T}(J^-X,J^-Y),\\
\mathrm{T}(JX,JY)&=& \mathrm{T}(J^+X,J^+Y)-\mathrm{T}(J^+X,J^-Y)-\mathrm{T}(J^-X,J^+Y)+\mathrm{T}(J^-X,J^-Y)\\
                                                 &=& \mathrm{T}(J^+X,J^+Y)+\mathrm{T}(J^-X,J^-Y),
\end{eqnarray*}
thus proving $\mathrm{T}(JX,JY)= \mathrm{T}(X,Y)$. $\blacksquare$
\bigskip

In \cite[Prop. 3.1]{etayo} Cruceanu and one of us had defined a  connection in an almost para-Hermitian manifold as the unique natural connection satisfying condition $ii)$ in the above lemma. Then, we can define the Chern-type connection on a $(J^2=\pm1)$-metric manifold satisfying $\alpha \varepsilon=-1$ as the connection determined in the following:

\begin{teor}[{\cite[Prop. 3.1]{etayo}, \cite[Theor. 6.1]{gray}}]
\label{teor:chern-connection}
Let  $(M,J,g)$ be a  $(J^2=\pm1)$-metric manifold with  $\alpha\varepsilon =-1$. Then there exist a unique linear connection $\nabla^{\mathrm{c}}$ in $M$ reducible to the $G_{(\alpha,\varepsilon)}$-structure defined by $(J,g)$  whose torsion tensor $\mathrm{T}^{\mathrm{c}}$ satisfies
\[
\mathrm{T}^{\mathrm{c}}(JX,JY)= \alpha \mathrm{T}^{\mathrm{c}}(X,Y), \quad \forall X, Y \in \mathfrak{X} (M).
\]
This connection will be named the Chern connection of the manifold $(M,g,J)$.
\end{teor}

\begin{obs}
\label{noChern}
We are going to check that in a $(1,1)$-metric manifold there is no a unique reducible connection satisfying
\[
 \mathrm{T}(X,Y) =0, \quad \forall X \in T^+_J (M), \forall Y \in T^-_J(M).
\]
Observe that Lemma \ref{teor:chernparahermitian} is also true in the case of an $(1,1)$-metric manifold, because the metric does not appear in the result. Thus, this Remark shows that a Chern connection can not be defined in  a $(1,1)$-metric manifold.
\bigskip

Let us prove the result. Taking an adapted local frame $(X_1, \ldots, X_n, Y_1,\ldots, Y_n)$ to the $G_{(1,1)}$-structure in $U$, one has
 \[
\mathrm{T}(X_i,Y_j)= \nabla_{X_i} Y_j -\nabla_{Y_j} X_i -[X_i, Y_j]=0, \quad \forall i, j =1, \ldots, n.
\]
As the linear connection $\nabla$ is reducible, then it is determined in  $U$ by the following functions
\[
 \nabla_{X_i} X_j =\sum_{k=1}^n \Gamma_{ij}^k X_k, \quad  \nabla_{Y_i} X_j = \sum_{k=1}^n \bar \Gamma_{ij}^k X_k, \quad
  \nabla_{X_i} Y_j =\sum_{k=1}^n \Gamma_{ij}^{k+n} Y_k, \quad   \nabla_{Y_i} Y_j =\sum_{k=1}^n \bar  \Gamma_{ij}^{k+n} Y_k, \quad i, j, k =1, \ldots, n.
\]
thus obtaining
\begin{eqnarray*}
g(\nabla_{X_i} Y_j, X_k) - g(\nabla_{Y_j} X_i, X_k) &=& g([X_i, Y_j], X_k),\\
g(\nabla_{X_i} Y_j, Y_k) - g(\nabla_{Y_j} X_i, Y_k) &=& g([X_i, Y_j], Y_k), \quad \forall i, j, k =1, \ldots, n,
\end{eqnarray*}
and then
\[
\Gamma_{ij}^{k+n} = g([X_i,Y_j],Y_k),\quad \bar \Gamma_{ij}^k = g([Y_j,X_i],X_k),  \quad \forall i, j, k =1, \ldots, n.
\]

These equalities do not impose any condition on the functions
\[
\Gamma_{ij}^k, \bar \Gamma_{ij}^{k+n}, \quad i,j,k=1, \ldots, n,
\]
and then, the condition is not enough to determine a unique linear connection.
\end{obs}

Finally, we will obtain the relation between the Chern connection  $\nabla^{\mathrm{c}}$ and the well adapted connection  $\nabla^{\mathrm{w}}$, in the case $\alpha \varepsilon =-1$. We will use  the potential tensor of
 $\nabla^{\mathrm{w}}$ , which is given by
 \[
S^{\mathrm{w}}(X,Y)=\nabla^{\mathrm{w}}_X Y -\nabla^{\mathrm{g}}_X Y, \quad \forall X, Y \in \mathfrak{X} (M).
\]
according to Definition \ref{teor:tensorpotencial}. The first result we need is the following:

\begin{prop}
Let $(M,J,g)$ be a  $(J^2=\pm1)$-metric manifold. The following conditions are equivalent:

\begin{enumerate}
\renewcommand*{\theenumi}{\roman{enumi})}
\renewcommand*{\labelenumi}{\theenumi}

\item $\mathrm{T}^{\mathrm{w}}(X,Y) +\varepsilon \mathrm{T}^{\mathrm{w}}(JX,JY)=0$, for all vector fields $X,Y$ in $M$.

\item $S^{\mathrm{w}}(X,Y)= \frac{(-\alpha)}{2} (\nabla^{\mathrm{g}}_X J) JY$,  for all vector fields $X,Y$ in $M$.

\end{enumerate}

\end{prop}

{\bf Proof.} As $\nabla^{\mathrm{w}}$ is a natural connection, by Lemma \ref{teor:natural}, one has
\[
\nabla^{\mathrm{w}}_X Y = \nabla^{\mathrm{g}}_X Y + S^{\mathrm{w}} (X,Y),  \quad \forall X,Y \in \mathfrak{X} (M).
\]

Taking into account $\nabla^{\mathrm{g}}$ is torsionless, one easily checks the following relation

\[
 \mathrm{T}^{\mathrm{w}}(X,Y)=S^{\mathrm{w}}(X,Y)-S^{\mathrm{w}}(Y,X), \quad \forall X,Y \in \mathfrak{X} (M).
\]
Given vector fields $X,Y,Z$ in $M$ one has
\begin{eqnarray*}
g(\mathrm{T}^{\mathrm{w}}(X,Y),Z)-g(\mathrm{T}^{\mathrm{w}}(Z,Y),X)&=& g(S^{\mathrm{w}}(X,Y),Z)-g(S^{\mathrm{w}}(Y,X),Z)\\
                                                 &-&g(S^{\mathrm{w}}(Z,Y),X)+g(S^{\mathrm{w}}(Y,Z),X)\\
                                                 &=& -g(S^{\mathrm{w}}(X,Z),Y)+g(S^{\mathrm{w}}(Z,X),Y)\\
                                                 &+&g(S^{\mathrm{w}}(Y,Z),X)+g(S^{\mathrm{w}}(Y,Z),X)\\
                                                 &=& g(\mathrm{T}^{\mathrm{w}}(Z,X),Y)+2g(S^{\mathrm{w}}(Y,Z),X),\\
g(\mathrm{T}^{\mathrm{w}}(JX,Y),JZ)-g(\mathrm{T}^{\mathrm{w}}(JZ,Y),JX)&=& g(\mathrm{T}^{\mathrm{w}}(JZ,JX),Y)+2g(S^{\mathrm{w}}(Y,JZ),JX),\\
                                                        &=& g(\mathrm{T}^{\mathrm{w}}(JZ,JX),Y)+2\alpha\varepsilon g(JS^{\mathrm{w}}(Y,JZ),X),
\end{eqnarray*}
Taking into account the above formulas and  formula (\ref{eq:welladapted}) one obtains the following relation satisfied by the well adapted connection:
\[
g(\mathrm{T}^{\mathrm{w}}(Z,X)+\varepsilon\mathrm{T}^{\mathrm{w}}(JZ,JX),Y)=-2g(S^{\mathrm{w}}(Y,Z)+\alpha JS^{\mathrm{w}}(Y,JZ),X), \quad \forall X,Y, Z \in\mathfrak{X} (M).
\]
Then one has the following chain of equivalences:
\begin{eqnarray*}
\mathrm{T}^{\mathrm{w}}(Z,X)+\varepsilon \mathrm{T}^{\mathrm{w}}(JZ,JX) =0 &\Leftrightarrow& S^{\mathrm{w}}(Y,Z)+\alpha JS^{\mathrm{w}}(Y,JZ)=0\\
                                                               &\Leftrightarrow& JS^{\mathrm{w}}(Y,Z)+S^{\mathrm{w}}(Y,JZ)=0\\
                                                               &\Leftrightarrow& JS^{\mathrm{w}}(Y,Z)=- S^{\mathrm{w}}(Y,JZ)\\
                                                               &\Leftrightarrow& (\nabla^{\mathrm{g}}_Y J)Z = S^{\mathrm{w}}(Y,Z)-S^{\mathrm{w}}(Y,JZ)=2JS^{\mathrm{w}}(Y,Z)\\
                                                               &\Leftrightarrow& S^{\mathrm{w}}(Y,Z)=\frac{(-\alpha)}{2} (\nabla^{\mathrm{g}}_Y J) JZ,
\end{eqnarray*}
for all vector fields  $X,Y,Z$ in $M$. Then one has
\begin{equation}
\mathrm{T}^{\mathrm{w}}(X,Y) +\varepsilon \mathrm{T}^{\mathrm{w}}(JX,JY)=0 \Leftrightarrow \nabla^{\mathrm{w}}_X Y = \nabla^{\mathrm{g}}_X Y + \frac{(-\alpha)}{2} (\nabla^{\mathrm{g}}_X J) JY, \quad \forall X,Y \in \mathfrak{X} (M). \  \blacksquare
\label{eq:nw=no}
\end{equation}

\begin{teor}
\label{teor:quasi}
Let $(M,J,g)$ be a $(J^2=\pm1)$-metric manifold with  $\alpha\varepsilon=-1$. Then the well adapted connection and the Chern connection  coincide if and only if
\[
\nabla^{\mathrm{w}}_X Y = \nabla^{\mathrm{g}}_X Y + \frac{(-\alpha)}{2} (\nabla^{\mathrm{g}}_X J) JY, \quad \forall X,Y \in \mathfrak{X} (M).
\]
\end{teor}

{\bf Proof. } As $\alpha ,\varepsilon \in \{ -1,1\} $ and $\alpha\varepsilon=-1$, then $\alpha =-\varepsilon $, and then equation (\ref{eq:nw=no}) reads as
\[
\mathrm{T}^{\mathrm{w}}(X,Y) -\alpha \mathrm{T}^{\mathrm{w}}(JX,JY)=0 \Leftrightarrow \nabla^{\mathrm{w}}_X Y = \nabla^{\mathrm{g}}_X Y + \frac{(-\alpha)}{2} (\nabla^{\mathrm{g}}_X J) JY, \quad \forall X,Y \in \mathfrak{X} (M),
\]
and the result trivially follows from Theorem \ref{teor:chern-connection}. $\blacksquare$
\bigskip

If $(M,J,g)$ is a $(J^2=\pm1)$-metric manifold with $\alpha\varepsilon=-1$, the connection given by the law derivation
\[
\nabla^{0}_X Y = \nabla^{\mathrm{g}}_X Y + \frac{(-\alpha)}{2} (\nabla^{\mathrm{g}}_X J) JY, \quad \forall X,Y \in \mathfrak{X} (M),
\]
is called the first canonical connection of $(M,J,g)$. Then the Chern connection and the well adapted connection are the same if and only if both connections coincide with the first canonical connection of $(M,J,g)$.

\end{document}